\def\date{22.5.07}

\input amssym.def
\input amssym.tex

\def\item#1{\vskip1.3pt\hang\textindent {\rm #1}}


\newskip\litemindent
\litemindent=0.7cm  
\def\Litem#1#2{\par\noindent\hangindent#1\litemindent
\hbox to #1\litemindent{\hfill\hbox to \litemindent
{\ninerm #2 \hfill}}\ignorespaces}
\def\litem{\Litem1}

\tolerance=300
\pretolerance=200
\hfuzz=1pt
\vfuzz=1pt

\hoffset=0in
\voffset=0.5in

\hsize=5.8 true in 
\vsize=9.2 true in
\parindent=25pt
\mathsurround=1pt
\parskip=1pt plus .25pt minus .25pt
\normallineskiplimit=.99pt

\countdef\revised=100
\mathchardef\emptyset="001F 
\chardef\ss="19
\def\3{\ss}
\def\anf{$\lower1.2ex\hbox{"}$}
\def\frac#1#2{{#1 \over #2}}
\def\>{>\!\!>}
\def\<{<\!\!<}

\def\into{\hookrightarrow}
\def\onto{\to\mskip-14mu\to} 
\def\ssssarr{\hbox to 15pt{\rightarrowfill}}
\def\sssarr{\hbox to 20pt{\rightarrowfill}}
\def\ssarr{\hbox to 30pt{\rightarrowfill}}
\def\sarr{\hbox to 40pt{\rightarrowfill}}
\def\arr{\hbox to 60pt{\rightarrowfill}}
\def\larr{\hbox to 60pt{\leftarrowfill}}
\def\Arr{\hbox to 80pt{\rightarrowfill}}

\def\ssssmapright#1{\smash{\mathop{\ssssarr}\limits^{#1}}}
\def\sssmapright#1{\smash{\mathop{\sssarr}\limits^{#1}}}
\def\ssmapright#1{\smash{\mathop{\ssarr}\limits^{#1}}}

\def\Alt{\mathop{\rm Alt}\nolimits}
\def\ad{\mathop{\rm ad}\nolimits}

\def\Ad{\mathop{\rm Ad}\nolimits}

\def\Alt{\mathop{\rm Alt}\nolimits}

\def\Der{\mathop{\rm Der}\nolimits}
\def\deg{\mathop{\rm deg}\nolimits}

\def\det{\mathop{\rm det}\nolimits}

\def\Diff{\mathop{\rm Diff}\nolimits}

\def\Exp{\mathop{\rm Exp}\nolimits}

\def\ev{\mathop{\rm ev}\nolimits}

\def\End{\mathop{\rm End}\nolimits}

\def\GL{\mathop{\rm GL}\nolimits}

\def\Hom{\mathop{\rm Hom}\nolimits}%
\def\id{\mathop{\rm id}\nolimits} 
\def\im{\mathop{\rm im}\nolimits}


%


\def\sgn{\mathop{\rm sgn}\nolimits}

\def\span{\mathop{\rm span}\nolimits}

\def\Sym{\mathop{\rm Sym}\nolimits}


\def\Tr{\mathop{\rm Tr}\nolimits}
\def\trile{\trianglelefteq}

\def\0{{\bf 0}}
\def\1{{\bf 1}}

\def\d{{\frak d}}

\def\g{{\frak g}}
\def\gl{{\frak {gl}}}
\def\h{{\frak h}}

\def\k{{\frak k}}

\def\n{{\frak n}}

\def\s{{\frak s}}

\def\z{{\frak z}}

\def\L{\mathop{\bf L{}}\nolimits}

\def\C{{{\Bbb C}{\mskip+1mu}}} 

\def\R{{\Bbb R}} 
\def\Z{{\Bbb Z}} 
\def\N{{\Bbb N}}

\def\F{{\Bbb F}}

\def\SS{{\Bbb S}} 
\def\T{{\Bbb T}} 

\def\:{\colon}  
\def\.{{\cdot}}
\def\|{\Vert}
\def\bsk{\bigskip}

\def\giantskip{\vskip2\bigskipamount}
\def\gsk{\giantskip}

\def\msk{\medskip}

\def \res {\!\mid\!\!}

\def\bbr{\bigbreak}
\def\giantbreak{\par \ifdim\lastskip<2\bigskipamount \removelastskip
         \penalty-400 \giantskip\fi}

\def\nin{\noindent}
\def\cen{\centerline}
\def\pagebreak{\vskip 0pt plus 0.0001fil\break}
\def\linebreak{\break}

\def\hat{\widehat}

\def\derat#1{{d \over dt} \hbox{\vrule width0.5pt 
                height 5mm depth 3mm${{}\atop{{}\atop{\scriptstyle t=#1}}}$}}

\def\epsilon{\varepsilon}

\def\nin{\noindent}
\def\oline{\overline}

\def\pder#1,#2,#3 { {\partial #1 \over \partial #2}(#3)}
\def\pde#1,#2 { {\partial #1 \over \partial #2}}
\def\phi{\varphi}


\def\subeq{\subseteq}

\def\tilde{\widetilde}

\font\ninerm=cmr9
\font\eightrm=cmr8

\font\eightbf=cmbx8


\font\smc=cmcsc10
\font\bfone=cmbx10 scaled\magstep1 
\font\bftwo=cmbx10 scaled\magstep2 

\def\qed{{\unskip\nobreak\hfil\penalty50\hskip .001pt \hbox{}\nobreak\hfil
          \vrule height 1.2ex width 1.1ex depth -.1ex
           \parfillskip=0pt\finalhyphendemerits=0\medbreak}\rm}

\def\qeddis{\eqno{\vrule height 1.2ex width 1.1ex depth -.1ex} $$
                   \medbreak\rm}

\def\Lemma #1. {\bigbreak\vskip-\parskip\noindent{\bf Lemma #1.}\quad\it}

\def\Sublemma #1. {\bigbreak\vskip-\parskip\noindent{\bf Sublemma #1.}\quad\it}

\def\Proposition #1. {\bigbreak\vskip-\parskip\noindent{\bf Proposition #1.}
\quad\it}

\def\Corollary #1. {\bigbreak\vskip-\parskip\nin{\bf Corollary #1.}
\quad\it}

\def\Theorem #1. {\bigbreak\vskip-\parskip\noindent{\bf Theorem #1.}
\quad\it}

\def\Definition #1. {\rm\bigbreak\vskip-\parskip\noindent
{\bf Definition #1.}
\quad}

\def\Remark #1. {\rm\bigbreak\vskip-\parskip\noindent{\bf Remark #1.}\quad}

\def\Example #1. {\rm\bigbreak\vskip-\parskip\noindent{\bf Example #1.}\quad}
\def\Examples #1. {\rm\bigbreak\vskip-\parskip\noindent{\bf Examples #1.}\quad}

\def\Problems #1. {\bigbreak\vskip-\parskip\noindent{\bf Problems #1.}\quad}
\def\Problem #1. {\bigbreak\vskip-\parskip\noindent{\bf Problem #1.}\quad}
\def\Exercise #1. {\bigbreak\vskip-\parskip\noindent{\bf Exercise #1.}\quad}

\def\Conjecture #1. {\bigbreak\vskip-\parskip\noindent{\bf Conjecture #1.}\quad}

\def\Proof#1.{\rm\par\ifdim\lastskip<\bigskipamount\removelastskip\fi\smallskip
            \noindent {\bf Proof.}\quad}

\def\Axiom #1. {\bigbreak\vskip-\parskip\noindent{\bf Axiom #1.}\quad\it}

\def\Satz #1. {\bigbreak\vskip-\parskip\noindent{\bf Satz #1.}\quad\it}

\def\Korollar #1. {\bbr\vskip-\parskip\nin{\bf Korollar #1.} \quad\it}

\def\Folgerung #1. {\bbr\vskip-\parskip\nin{\bf Folgerung #1.} \quad\it}

\def\Folgerungen #1. {\bbr\vskip-\parskip\nin{\bf Folgerungen #1.} \quad\it}

\def\Bemerkung #1. {\rm\bigbreak\vskip-\parskip\noindent{\bf Bemerkung #1.}
\quad}

\def\Beispiel #1. {\rm\bigbreak\vskip-\parskip\noindent{\bf Beispiel #1.}\quad}
\def\Beispiele #1. {\rm\bigbreak\vskip-\parskip\noindent{\bf Beispiele #1.}\quad}
\def\Aufgabe #1. {\rm\bigbreak\vskip-\parskip\noindent{\bf Aufgabe #1.}\quad}
\def\Aufgaben #1. {\rm\bigbreak\vskip-\parskip\noindent{\bf Aufgabe #1.}\quad}

\def\Beweis#1. {\rm\par\ifdim\lastskip<\bigskipamount\removelastskip\fi
           \smallskip\noindent {\bf Beweis.}\quad}

\nopagenumbers

\def\date{\ifcase\month\or January\or February \or March\or April\or May
\or June\or July\or August\or September\or October\or November
\or December\fi\space\number\day, \number\year}

\def\title{Title ??}
\def\author{Author ??}

\def\thanks#1{\footnote*{\eightrm#1}}

\def\rightheadline{\hfil{\eightrm\title}\hfil\tenbf\folio}
\def\leftheadline{\tenbf\folio\hfil{\eightrm\author}\hfil}
\headline={\vbox{\line{\ifodd\pageno\rightheadline\else\leftheadline\fi}}}

\def\firstheadline{}
\def\firstfootline{\cen{\rm\folio}}

\def\seite #1 {\pageno #1
               \headline={\ifnum\pageno=#1 \firstheadline
               \else\ifodd\pageno\rightheadline\else\leftheadline\fi\fi}
               \footline={\ifnum\pageno=#1 \firstfootline\else{}\fi}}

\newdimen\dimenone
 \def\checkleftspace#1#2#3#4{
 \dimenone=\pagetotal
 \advance\dimenone by -\pageshrink   
 \ifdim\dimenone>\pagegoal          
   \else\dimenone=\pagetotal
        \advance\dimenone by \pagestretch
        \ifdim\dimenone<\pagegoal
          \dimenone=\pagetotal
          \advance\dimenone by#1         
          \setbox0=\vbox{#2\parskip=0pt                
                     \hyphenpenalty=10000
                     \rightskip=0pt plus 5em
                     \noindent#3 \vskip#4}    
        \advance\dimenone by\ht0
        \advance\dimenone by 3\baselineskip   
        \ifdim\dimenone>\pagegoal\vfill\eject\fi
          \else\eject\fi\fi}


\def\subheadline #1{\nin\bigbreak\vskip-\lastskip
      \checkleftspace{0.9cm}{\bf}{#1}{\medskipamount}
          \indent\vskip0.7cm\centerline{\bf #1}\medskip}
\def\subsection{\subheadline} 

\def\lsubheadline #1 #2{\bigbreak\vskip-\lastskip
      \checkleftspace{0.9cm}{\bf}{#1}{\bigskipamount}
         \vbox{\vskip0.7cm}\cen{\bf #1}\msk \cen{\bf #2}\bsk}

\def\sectionheadline #1{\bigbreak\vskip-\lastskip
      \checkleftspace{1.1cm}{\bf}{#1}{\bigskipamount}
         \vbox{\vskip1.1cm}\cen{\bfone #1}\bsk}
\def\section{\sectionheadline} 

\def\lsectionheadline #1 #2{\bigbreak\vskip-\lastskip
      \checkleftspace{1.1cm}{\bf}{#1}{\bigskipamount}
         \vbox{\vskip1.1cm}\cen{\bfone #1}\msk \cen{\bfone #2}\bsk}

\def\lchapterheadline #1 #2{\bigbreak\vskip-\lastskip\indent\vskip3cm
                       \cen{\bftwo #1} \msk \cen{\bftwo #2} \gsk}
\def\llsectionheadline #1 #2 #3{\bigbreak\vskip-\lastskip\indent\vskip1.8cm
\cen{\bfone #1} \msk \cen{\bfone #2} \msk \cen{\bfone #3} \nobreak\bsk\nobreak}


\newtoks\literat
\def\[#1 #2\par{\literat={#2\unskip.}%
\hbox{\vtop{\hsize=.15\hsize\nin [#1]\hfill}
\vtop{\hsize=.82\hsize\nin\the\literat}}\par
\vskip.3\baselineskip}

\def\references{
\sectionheadline{\bf References}
\frenchspacing

\entries\par}

\mathchardef\emptyset="001F 
\def\address{Author: \tt$\backslash$def$\backslash$address$\{$??$\}$}

\def\abstract #1{{\narrower\baselineskip=10pt{\noindent
\eightbf Abstract.\quad \eightrm #1 }
\bigskip}}

\def\firstpage{\nin
{\obeylines \parindent 0pt }
\vskip2cm
\centerline{\bfone\title}
\gsk
\centerline{\bf\author}
\vskip1.5cm \rm}

\def\addresstwo{}

\def\dlastpage{\par\vbox{\vskip1cm\nin
\line{
\vtop{\hsize=.5\hsize{\parindent=0pt\baselineskip=10pt\nin\address}}
\quad 
\vtop{\hsize=.42\hsize\nin{\parindent=0pt
\baselineskip=10pt\addresstwo}}
\hfill} }}

\def\Box #1 { \msk\par\nin 
\centerline{
\vbox{\offinterlineskip
\hrule
\hbox{\vrule\strut\hskip1ex\hfil{\smc#1}\hfill\hskip1ex}
\hrule}\vrule}\msk }

\def\adots{\mathinner{\mkern1mu\raise1pt\vbox{\kern7pt\hbox{.}}
                        \mkern2mu\raise4pt\hbox{.}
                        \mkern2mu\raise7pt\hbox{.}\mkern1mu}}


\pageno=1
\def\title{On the cohomology of vector fields on parallelizable manifolds} 
\def\author{Yuly Billig, Karl-Hermann Neeb} 


\def\div{\mathop{\rm div}\nolimits}%

\def\half{\textstyle{1\over 2}}

\def\V{{\cal V}}
\def\WW{{\widehat W}_N}
\def\tW{{\widetilde W}(N)}
\def\L{{\cal L}}
\def\Lie{\mathop{\bf L{}}\nolimits}

\def\Fl{{\rm Fl}}
\def\equ{{\rm eq}}

\def\bp{{\overline \psi}}
\def\Bp{{\overline \Psi}}

\def\d{\partial}
\def\F{{\cal F}}

\def\fd{{\frak d}}
\def\sst{\scriptstyle} 

\firstpage 

\abstract{In the present paper we determine for each parallelizable 
smooth compact manifold ${\sst M}$ the cohomology spaces 
${\sst H^2({\cal V}_M,\oline\Omega^p_M)}$ of the Lie algebra 
${\sst {\cal V}_M}$ of smooth vector fields on ${\sst M}$ with values in the 
module ${\sst \oline\Omega^p_M = \Omega^p_M/d\Omega^{p-1}_M}$. 
The case of ${\sst p=1}$ is of particular interest since the gauge algebra
${\sst C^\infty (M,\k)}$ has the universal central extension with center
${\sst \oline\Omega^1_M}$, 
generalizing affine Kac-Moody algebras. The second
cohomology ${\sst H^2(\V_M, \oline\Omega^1_M)}$ classifies twists of the
semidirect product of ${\sst \V_M}$ with 
the universal central extension ${\sst C^\infty (M,\k) \oplus 
\oline\Omega^1_M}$.
\hfill\break 
Keywords: Lie algebra of vector fields, Lie algebra cohomology, 
Gelfand-Fuks cohomology, extended affine Lie algebra \hfill\break 
}

One of the most important insights in the theory of affine 
Kac--Moody Lie algebras is that they can all be realized as 
twisted or untwisted loop algebras. In the untwisted case, this 
realization starts with the Lie algebra of maps from a circle $\SS^1$
to a finite-dimensional complex simple Lie algebra $\k$, 
$$ {\cal L}(\k) := \C[t,t^{-1}] \otimes \k, $$
then proceeds with the universal central extension 
$$ \C \into \hat{\cal L}(\k) \onto {\cal L}(\k), $$
and completes the picture by adding a derivation $d$ acting by 
$t {d\over dt}$ on ${\cal L}(\k)$ and, accordingly, on the 
central extension. In the representation theory of affine algebras, 
an important role is played by the Virasoro algebra, which emerges 
via the Sugawara construction. At the Lie algebra level,
the Witt algebra of vector fields on a circle,
$\fd := \Der(\C[t,t^{-1}])$ acts on $\hat{\cal L}(\k)$, so that we may 
form the semidirect product 
$\g := \hat{\cal L}(\k) \rtimes \fd.$
To get the Virasoro algebra, we need to twist the above semidirect product
by a $2$-cocycle $\tau \in Z^2(\fd,\C)$, which leads to the affine-Virasoro Lie algebra $\g_\tau$ 
with the bracket 
$$ [(x,d),(x',d')] = ([x,x'] + d.x' - d'.x + \tau(d,d'), [d,d']), $$
which contains both affine Kac-Moody algebra and the Virasoro algebra as
subalgebras.

The theory of affine Kac--Moody algebras has an analytic side, 
where one replaces the algebra $\C[t,t^{-1}]$ of Laurent polynomials 
by the Fr\'echet algebra $C^\infty(\SS^1,\C)$ of complex valued 
functions on the circle. In this context one also obtains 
a one-dimensional central extension and the role of the Witt algebra 
is played by the Fr\'echet--Lie algebra ${\cal V}_{\SS^1}$ of smooth 
vector fields on the circle, which also has a one-dimensional 
central extension. For an exposition of these ideas we refer to the monograph 
of Pressley and Segal [PS86].

 Let us discuss the generalization of this construction from $\SS^1$ to the case 
of an arbitrary $C^\infty$-manifold $M$. The gauge Lie algebra of $\k$-valued functions on $M$,
$$ \L_M(\k) = C^\infty (M,\k),$$
endowed with the pointwise 
defined Lie bracket and the natural Fr\'echet topology,
has the universal central extension with the central space
$$\oline\Omega^1(M,\C) := \Omega^1(M,\C)/d\Omega^0(M,\C)$$
(cf.\ [Ka84], [Ma02]). The Lie bracket on the Lie algebra
$$ \hat\L_M(\k) = C^\infty (M,\k) \oplus \oline\Omega^1(M,\C) $$
is given by the formula
$$[f_1 \otimes g_1, f_2 \otimes g_2] = f_1 f_2 \otimes [g_1,g_2]
+ (g_1 | g_2) [f_2 d f_1],$$ 
where $f_1, f_2 \in C^\infty (M, \C)$, $g_1, g_2 \in \k$, $( \cdot | \cdot)$
is the Cartan--Killing form on $\k$ and $[\alpha]$ denotes the class 
of a $1$-form $\alpha$ in $\oline\Omega^1(M,\C)$.

 The Lie algebra ${\cal V}_M$ of smooth vector fields on $M$ acts on both
the gauge algebra $C^\infty (M, \k)$ and on $\oline\Omega^1(M,\C)$ by the 
Lie derivative. 
This action is compatible with the above central cocycle,
allowing us to consider the semidirect product of ${\cal V}_M$ with $ \hat\L_M(\k)$:
$$\g = \Big( C^\infty (M,\k) \oplus \oline\Omega^1(M,\C) \Big)
\rtimes {\cal V}_M .$$
However the interplay between affine Lie algebras and the Virasoro algebra suggests that 
it will be natural here to twist the Lie bracket in $\g$ by means of a 2-cocycle
$\tau \in 
Z^2_c({\cal V}_M,\oline\Omega^1(M,\C))$, 
resulting in the family of Lie algebras
$$\g_\tau = \Big(C^\infty (M,\k) \oplus \oline\Omega^1(M,\C)\Big) 
\oplus_\tau {\cal V}_M ,$$
so that the equivalence 
classes of twistings are classified by  the cohomology space 
$$ H^2_c({\cal V}_M,\oline\Omega^1(M,\C)) 
\cong H^2_c({\cal V}_M,\oline\Omega^1(M,\R)) \otimes_\R \C. $$

In the present paper we calculate this space for 
parallelizable manifolds $M$, relying heavily on results on 
the cohomology of Lie algebras of vector fields with values 
in differential forms by Gelfand--Fuks, Haefliger and Tsujishita. 
We calculate this space using the short exact sequences 
$$ \0 \to H^1_{\rm dR}(M,\R) \into \oline\Omega^1(M,\R)) \to B^2_{\rm dR}(M,\R) \to \0 $$
and 
$$ \0 \to B^1_{\rm dR}(M,\R) \into \Omega^1(M,\R) \to \oline\Omega^1(M,\R) \to \0, 
$$
and then applying the corresponding long exact sequence in cohomology. 
Hence we need detailed knowledge on the cohomology of 
${\cal V}_M$ with values in differential forms. 

\bsk 

In the following $M$ denotes an 
$N$-dimensional parallelizable compact manifold.  
The main goal of this paper is to describe for each $p \in \N$ 
the second cohomology spaces 
$$ H^2_c({\cal V}_M, \oline\Omega^p(M,\R)) 
\quad \hbox{ for}\quad  p \in \N, \quad 
\oline\Omega^p(M,\R) := \Omega^p(M,\R)/d\Omega^{p-1}(M,\R) $$
with respect to the natural action of the Lie algebra ${\cal V}_M$ of smooth 
vector fields. For reasons described above, our main interest lies 
in the case $p = 1$. For $N = 1$ and connected $M$ 
we have $M \cong \SS^1$ and $\oline\Omega^1(M,\R) \cong \R$ is 
a trivial module. Then $H^2_c({\cal V}_{\SS^1},\R) \cong \R$ describes 
the central extensions of ${\cal V}_{\SS^1}$, and a generator 
corresponds to the Virasoro algebra. 
For $N \geq 2$, we find that 
$$ H^2_c({\cal V}_M, \oline\Omega^1(M,\R)) 
\cong H^3_{\rm dR}(M,\R) \oplus \R^2. $$
Here we use that each closed $3$-form $\omega$ yields 
an $\oline\Omega^1(M,\R)$-valued $2$-cocycle by 
$\omega^{[2]}(X,Y) := [i_Y i_X\omega]$ (cf.\ [Ne06a]). In terms of a 
trivializing $1$-form $\kappa \in \Omega^1(M,\R^N)$, the other 
two cocycles can be described as follows. Define 
$\theta \: {\cal V}_M \to C^\infty(M,\gl_N(\R))$ by 
${\cal L}_X\kappa = - \theta(X)\cdot \kappa$, 
$$ \oline\Psi_1(X) := \Tr(\theta(X)) \in \Omega^0(M,\R) = C^\infty(M,\R), \quad 
\Psi_1(X) := \Tr(d\theta(X)) \in \Omega^1(M,\R),$$
and 
$$ \oline\Psi_2(X,Y) := 
[\Tr(d\theta(X) \wedge \theta(Y) - d\theta(Y) \wedge \theta(X))] 
\in \oline\Omega^1(M,\R). $$
Then $\oline\Psi_1 \wedge \Psi_1$ and $\oline\Psi_2$ provide the two 
additional generators of $H^2_c({\cal V}_M,\oline\Omega^1(M,\R))$. 

 In the case when $M$ is a torus $\T^N$, we give explicit formulas for these cocycles in 
coordinates. For the toroidal Lie algebras, the cocycles $\oline\Psi_2$ and 
$\oline\Psi_1 \wedge \Psi_1$ were discovered in the representation theory of these
algebras ([EM94], [L99], see also [BB99], [Bi06]). It turns out here
that it is easier to construct representations of $\g_\tau$ for a non-trivial cocycle
$\tau$, rather than for the semidirect product $\g$.

 For the toroidal Lie algebras it would be natural to work in the algebraic setting,
considering the algebra $A$ 
of Laurent polynomials in $N$ variables as the algebra
of functions on $\T^N$ (Fourier polynomials), 
instead of $C^\infty (\T^N, \C)$. 
Although the cocycles that we get here are well-defined in this algebraic setting,
our results do not solve the problem of describing the algebraic cohomology 
$H^2(\Der(A), \Omega^1(A)/dA)$. The reason for this is that our proof
is based on the results of Tsujishita and Haefliger, while the algebraic counterparts of these
have not been established.

 Toroidal Lie algebras are also closely related to the class of extended affine Lie algebras
([Neh04], [BN06]). According to a result of [ABFP05], 
most of the extended affine Lie algebras
are twisted versions of toroidal Lie algebras. However, instead of the Lie algebra of
vector fields ${\cal V}_{\T^N}$, in the theory of extended affine Lie algebras one must use the Lie
algebra of divergence zero vector fields ${\cal V}^{div}_{\T^N}$ (or a subalgebra).
Thus it would be also interesting to describe the space $H^2 ({\cal V}^{div}_{\T^N},\oline\Omega^1(\T^N,\R))$.
This question remains open. We point out here that the restriction of the 
cocycle $\oline\Psi_1 \wedge \Psi_1$ to 
${\cal V}^{div}_{\T^N}$ vanishes, while $\oline\Psi_2$ does not. 
Also we note that the cocycles corresponding to $H^3_{\rm dR}(M,\R)$ are not relevant for the
theory of extended affine Lie algebras either, because of the additional restriction that the cocycle should 
vanish on the subalgebra of degree zero derivations, so that this subalgebra remains abelian.

\msk 

The cohomology of the subcomplex $C^p_{\rm loc}({\cal V}_M,\Omega^p)$ 
of {\it local cochains}, i.e., cochains 
$f \in C^p_c({\cal V}_M,\Omega^p_M)$ that are differential operators in 
each argument is easier accessible than the full complex. 
In view of [Tsu81], the difference fully comes from the case $p = 0$. 
In this case, the description of 
$H^\bullet_{\rm loc}({\cal V}_M, {\cal F}_M)$ is facilitated significantly 
by the results of de Wilde and Lecomte ([dWL83]) who 
show that the cohomology of the differential graded algebra 
$C^\bullet_{\rm loc}({\cal V}_M, {\cal F}_M)$ coincides with the 
cohomology of the subalgebra generated by 
the image of $\Omega^\bullet_M$, the image of the Chern--Weil homomorphism 
$$ \chi_1 \: \Sym^\bullet(\gl_N(\R),\R)^{\gl_N(\R)} 
\to C^{2\bullet}_{\rm loc}({\cal V}(M), {\cal F}_M) $$
corresponding to a connection on the frame bundle $J^1(M)$ of $M$ 
(which yields cocycles depending on first order derivatives), and 
the image of an algebra homomorphism 
$$ \chi_2 \: C^\bullet(\gl_N(\R),\R) 
\to C^\bullet_{\rm loc}({\cal V}(M), {\cal F}_M) $$
whose image consists of cocycles if the connection is flat. 
If all Prontrjagin classes of $M$ vanish, then this results in an isomorphism 
$$H^\bullet_{\rm loc}({\cal V}_M, {\cal F}_M) 
\cong H^\bullet_{\rm dR}(M,\R) \otimes H^\bullet(\gl_N(\R),\R)
\cong H^\bullet_{\rm dR}(M,\R) \otimes C^\bullet(\gl_N(\R),\R)^{\gl_N(\R)}.$$
As our results show, for $N > 1$, all classes in 
$H^2_c({\cal V}_M, \oline\Omega^p(M,\R))$ can be represented by 
local cochains, but  the generator of 
$H^2_c({\cal V}_{\SS^1}, \R) \subeq H^2_c({\cal V}_{\SS^1}, C^\infty(M,\R))$ 
is non-local because it involves an integration. 


Some of Tsujishita's results have been generalized by Rosenfeld 
in [Ro71] to other classes of irreducible primitive Lie subalgebras 
of ${\cal V}_M$. It would be of some interest to see if these results 
could be used to obtain the cohomology of these Lie algebras with 
values in $\oline\Omega^1_M$.

\subheadline{Notation and conventions} 

If $V$ is a module of the Lie algebra $\g$, we write 
$C^p(\g,V)$ for the corresponding space $p$-cochains, 
$d_\g \: C^p(\g,V) \to C^{p+1}(\g,V)$ for the Chevalley--Eilenberg differential, 
$B^p(\g,V) = d_\g(C^{p-1}(\g,V))$ for the space of $p$-coboundaries, 
$Z^p(\g,V) = \ker(d_\g \: C^p(\g,V) \to C^{p+1}(\g,V))$ 
for the space of $p$-cocycles and $H^p(\g,V)$ for the Lie algebra cohomology. 
If $\g$ is a topological Lie algebra and $V$ a continuous $\g$-module, then 
$C^p_c(\g,V)$ etc.\ stands for the space of continuous cochains. 

For $x \in \g$ we have the insertion operator 
$i_x \: C^{p+1}(\g,V) \to C^p(\g,V)$ and there is a natural action of 
$\g$ on $C^p(\g,V)$, which is denoted by the operators ${\cal L}_x$, satisfying 
the Cartan relation
$$ {\cal L}_x = d_\g \circ i_x + i_x \circ d_\g. $$

For a subalgebra $\h \leq \g$ we write 
$$C^p(\g,\h,V) := \{ \alpha \in C^p(\g,V) \: (\forall x \in \h)\ i_x \alpha = 0, 
i_x(d_\g\alpha) = 0\} $$
for the relative cochains modulo $\h$, and accordingly 
$B^p(\g,\h,V)$, $Z^p(\g,\h,V)$ and $H^p(\g,\h,V)$. 

\msk If $G$ is a group, we denote the identity element by  $\1$, and for $g \in G$, we write 
\par\nin $\lambda_g \: G \to G, x \mapsto gx$ for the
{\it left multiplication} by $g$, 
\par\nin $\rho_g \: G \to G, x \mapsto xg$ for the {\it right multiplication} by 
$g$, 
\par\nin $m_G \: G \times G \to G, (x,y) \mapsto xy$ for 
the {\it multiplication map}, and 
\par\nin $\eta_G \: G \to G, x \mapsto x^{-1}$ for the {\it inversion}. 

\msk If $M$ is a smooth manifold, we write 
${\cal V}_M$ for the Lie algebra of smooth vector fields on $M$, 
${\cal F}_M := C^\infty(M,\R)$ for the algebra of smooth 
real-valued functions on $M$, $\Omega^k_M := \Omega^k(M,\R)$, the space of smooth 
real-valued $p$-forms on $M$, 
and 
$$ Z^k_M := \ker(d\res_{\Omega^k_M}), \quad 
B^k_M := d\Omega_M^{k-1}, \quad 
H^k_M := H^k_{\rm dR}(M,\R), \quad \hbox{ and } \quad \oline\Omega^k_M := \Omega^k_M/B^k_M. $$
We write $T^{(p,q)}(M) := T(M)^{\otimes p} \otimes T^*(M)^{\otimes q}$ for the 
tensor bundles over $M$ and  $\Gamma(T^{(p,q)}(M))$ for their spaces of smooth 
sections, i.e., the $(p,q)$-tensor fields. 

\sectionheadline{I. Crossed homomorphisms of Lie algebras and pull-backs} 

In this short first section, we collect some basic facts on crossed 
homomorphisms of Lie algebras which provide the tools used throughout the 
forthcoming sections. The main point is Theorem~I.7 on maps in 
Lie algebra cohomology defined by crossed homomorphisms. 

Let $\h$ and $\n$ be two Lie algebras, with $\h$ acting on $\n$ by derivations, i.e., we are
given a Lie algebra homomorphism
$\tau: \h \rightarrow \Der(\n).$

\Example I.1. (a) Our primary example would be 
$\h = {\cal V}_M$ and $\n = C^\infty(M,\g)$, 
where $\g$ is a finite-dimensional Lie algebra and the bracket on $\n$ is defined 
pointwise. Then 
$$ (\tau(v).f)(m) := (v.f)(m) := df(m)v(m) $$
defines a Lie algebra homomorphism from ${\cal V}_M$ to $\Der (C^\infty(M,\g))$.

(b) If $A$ is a commutative associative 
algebra and $\g$ a finite-dimensional Lie algebra, 
then $\n := A \otimes \g$ carries a Lie algebra structure defined by 
$$ [a \otimes x, a' \otimes x'] := aa' \otimes [x,x']. $$
Then we obtain an action of $\h := \Der(A)$ on $\g$ by 
$\tau(D).(a \otimes x) := (D.a) \otimes x.$ 
\qed

\Definition I.2. A linear map 
$\theta: \h \rightarrow \n$ is called {\it a crossed homomorphism}
if 
$$\theta([x,y]) = [\theta(x), \theta(y)] + \tau(x)\theta(y) - \tau(y)\theta(x)
\quad\quad {\rm for \ all \ } x,y\in \h.
\qeddis

\Remark I.3. (a) Note that $\theta$ is a crossed homomorphism if and only if the map 
$(\theta, \id_\h) \: \h \to \n \rtimes \h$ is a homomorphism of Lie algebras. 

In terms of the Lie algebra cohomology of $\h$
with values in the $\h$-module $\n$, the condition that 
$\theta$ is a crossed homomorphism can be written as 
$(d_\h \theta)(x,y) + [\theta(x),\theta(y)] = 0,$
i.e., $\theta$ satisfies the Maurer--Cartan equation 
$$ d_\h \theta + \half [\theta, \theta] = 0.$$

(b) If $\h \cong \n \rtimes_\eta \s$ is a semidirect product and 
the action of $\h$ on $\n$ is defined by $\tau(n,s) := \eta(s)$, then 
$$ [(n,s),(n',s')] = ([n,n'] + \eta(s).n' - \eta(s').n, [s,s']) $$
implies that $\theta \: \h \to \n, (n,s) \mapsto n$ is a crossed homomorphism. 

(c) The kernel of any crossed homomorphism is a subalgebra. 
\qed

\Remark I.4. Let $\rho \: \h \to \n$ be a homomorphism of Lie algebras. 
Then we obtain on $\n$ an action of $\h$ by 
$\tau(x).y := [\rho(x),y]$. Now a linear map $\theta \: \g \to \n$ is a crossed 
homomorphism if and only if $\rho' := \rho + \theta$ is a homomorphism of Lie algebras. 

This follows directly by comparing 
$$ \rho'([x,y]) = \rho([x,y]) + \theta([x,y]), $$
with 
$$ [\rho'(x), \rho'(y)] = [\rho(x), \rho(y)] + [\theta(x),\rho(y)] 
+ [\rho(x),\theta(y)] + [\theta(x),\theta(y)]. $$

It follows in particular that $\theta := -\rho$ is a crossed homomorphism. 
\qed

Let us show how a crossed homomorphism may be used to 
pull back cocycles on $\n$ to cocycles on $\h$.
First, we need to define the notion of an equivariant cochain.

\Definition I.5. Let $V$ be a module for the Lie algebra $\h$ and a trivial module for $\n$. A cochain $\phi \in C^k (\n,V)$ is called {\it equivariant} if 
$$x \phi (y_1, \ldots, y_k) = \sum_{j=1}^k \phi(y_1, \ldots, \tau(x)y_j, \ldots, y_k)
\quad\quad {\rm for \ all \ }\quad  x\in \h, y_1,\ldots, y_k\in \n, $$
which is equivalent to ${\cal L}_x \phi = 0$ for each $x \in \h$. 
\qed

Since $[{\cal L}_x, d_\n] = 0$ holds for each $x \in \h$ on $C^\bullet(\n,V)$, we have: 
\Lemma I.6. Equivariant 
cochains form a subcomplex $C^\bullet_{\equ}(\n,V)$ in $C^\bullet (\n,V)$.
\qed

By definition, the subcomplex of equivariant 
cochains yields the equivariant cohomology $H^\bullet_{\equ}(\n,V)$.
Identifying $\n$ with a subalgebra of the semidirect product $\n \rtimes \h$, 
we may identify $C^p(\n,V)$ with the subspace 
$$ \{ \alpha \in C^p(\n\rtimes \h,V)  \: (\forall x \in \h) \ i_x \alpha = 0\}. $$
Hence the relative cochain space 
$$C^p(\h \rtimes \n,\h,V) := \{ \alpha \in C^p(\n \rtimes \h,V) 
\: (\forall x \in \h)\ i_x \alpha = 0, 
i_x(d\alpha) = 0\} $$
can be identified with $C^p_{\rm eq}(\n,V)$ because 
if $i_x \alpha = 0$ holds for each $x \in \h$, then ${\cal L}_x \alpha = i_x d\alpha$, 
so that invariance is equivalent to $\alpha \in C^p(\n \rtimes \h, \h, V)$ 
(cf.\ [Fu86, p.~16]). 

\Theorem I.7. Let $V$ be an $\h$-module and consider it as a trivial $\n$-module. 
Further let $\theta$ be a crossed homomorphism of Lie algebras 
$\theta: \h \rightarrow \n$.
Then the map 
$$\theta^* : C^\bullet_\equ (\n,V) \rightarrow C^\bullet (\h,V), \quad 
\phi \mapsto \phi \circ (\theta \times \ldots \times \theta) $$
is a morphism of cochain complexes, hence induces a linear map 
of cohomology spaces \break $\theta^* : H^\bullet_\equ (\n,V) \rightarrow H^\bullet (\h,V)$.

\Proof. We have seen above that we may identify the complex 
$(C^\bullet_{\rm eq}(\n,V),d_\n)$ with the relative Lie algebra subcomplex 
$$ (C^\bullet(\n \rtimes \h,\h,V), d_{\n \rtimes \h}) 
\subeq (C^\bullet(\n \rtimes \h), d_{\n \rtimes\h}). $$
For $\phi \in C^p_{\rm eq}(\n,V)$, we write $\tilde\phi$ for the corresponding element 
of $C^p(\n \rtimes \h, V)$. Then we have 
$$ \theta^*\phi = (\theta,\id_\h)^*\tilde\phi, $$
and the assertion follows from the fact that $(\theta,\id_\h) \: \h \to \n \rtimes \h$ 
is a morphism of Lie algebras, 
hence induces a morphism of cochain complexes 
$C^\bullet(\n\rtimes \h, V) \to C^\bullet(\h,V).$
\qed

\subheadline{Some applications to the Lie algebra of vector fields} 

In this subsection $M$ denotes a finite-dimensional smooth manifold. 

\Proposition I.8. A $\g$-valued differential form $\theta\in \Omega^1(M,\g)$ defines an  
${\cal F}_M$-linear crossed homomorphism 
$$ {\cal V}_M \to C^\infty(M,\g), \quad X \mapsto i_X \theta $$
if and only if $\theta$ satisfies the Maurer--Cartan equation 
$d\theta + {1\over 2}[\theta,\theta] = 0.$

\Proof. Since the exterior differential on $\Omega^\bullet(M,\g)$ 
coincides with the ${\cal V}_M$-Lie algebra differential, this is an immediate 
consequence of Remark~I.3(a). 
\qed

Recall that the left Maurer--Cartan form $\kappa_G \in \Omega^1(G,\g)$ of a Lie group $G$ 
with Lie algebra $\g$ is defined by $\kappa_G(v) := T(\lambda_g^{-1}) v$ 
for $v \in T_g(G)$. 

\Corollary I.9. {\rm(cf.\ [BR73], Th.~3.1)} If $M$ is a smooth manifold, $\g$ a Lie algebra 
and $\kappa \in \Omega^1(M,\g)$ satisfies the Maurer--Cartan equation, then we have a 
map 
$$ \zeta \: C^\bullet(\g,\R) \to  C^\bullet_c({\cal V}_M, {\cal F}_M), 
\quad 
\zeta(\omega)(X_1,\ldots, X_p) := \omega(\kappa(X_1),\ldots, \kappa(X_p)), $$
whose range consists of ${\cal F}_M$-multilinear 
cocycles represented by differential forms. 
In particular, we obtain an algebra homomorphism 
$$ H^\bullet(\g,\R) \to H^\bullet_c({\cal V}_M, {\cal F}_M).$$

\Proof. In view of Theorem~I.7, it only remains to verify that $\zeta$ 
is compatible with the multiplication of cochains, but this is an immediate 
consequence of the definition. 
\qed 

\Corollary I.10. If $G$ is a Lie group with Lie algebra $\g$ and 
$\kappa_G \in \Omega^1(G,\g)$ the left Maurer--Cartan form, then 
$$ \theta \: {\cal V}_G  \to C^\infty(G,\g), \quad X \mapsto i_X \kappa_G $$
is a ${\cal F}_G$-linear crossed homomorphism which is bijective. 

\Proof. In view of Proposition~I.8, we 
only have to recall that $\kappa_G$ satisfies the Maurer--Cartan equation. 
\qed

\sectionheadline{II. Cocycles with values in differential forms} 

Let $M$ be a smooth paracompact $N$-dimensional manifold. 
In this section we first explain how affine connections can be 
used to define for each $k \in \N$ 
cocycles $\Psi_k\in Z^k_c({\cal V}_M, \Omega^k_M)$. 
If, in addition, the tangent bundle of $M$ is trivial, we obtain 
two more families of cocycles 
$\oline\Psi_k\in Z^k_c({\cal V}_M, \oline\Omega^{k-1}_M)$ 
and $\Phi_k\in Z_c^{2k-1} ({\cal V}_M, \F_M)$. These cocycles 
satisfy the relation $d \circ \oline\Psi_k = \Psi_k$ that will 
play a crucial role in our determination of the spaces 
$H^2_c({\cal V}_M, \oline\Omega^p_M)$ in Section IV.

The set of affine connections $\nabla$ on $M$ is an affine space whose 
tangent space is the space $\Gamma(T^{(1,2)}(M))$ of tensors of type $(2,1)$. 
The elements of this space are ${\cal F}_M$-bilinear maps 
${\cal V}_M \times {\cal V}_M \to {\cal V}_M$ and any such map 
can also be considered as an ${\cal F}_M$-linear map 
${\cal V}_M \to \Gamma(\End(TM))$, i.e., as a $1$-form with values in the 
endomorphism bundle $\End(TM) \cong T^{(1,1)}(M)$ of the tangent bundle $T(M)$. 
In this sense we identify 
$\Gamma(T^{(1,2)}(M))$ with the space $\Omega^1(M,\End(TM))$ of 
$1$-forms with values in the endomorphism bundle $\End(TM)$. 
Note that this space 
carries a natural module structure for the Lie algebra ${\cal V}_M$, 
given by the Lie derivative. 

\Lemma II.1. {\rm([Ko74])} Any affine connection $\nabla$ on $M$ defines a $1$-cocycle 
$$ \zeta \: {\cal V}_M \to \Omega^1(M,\End(TM)), \quad X \mapsto {\cal L}_X\nabla, $$
where 
$$ ({\cal L}_X\nabla)(Y)(Z) := [X,\nabla_Y Z] - \nabla_{[X,Y]}Z - \nabla_Y[X,Z]. $$
For any other affine connection $\nabla'$ the corresponding cocycle $\zeta'$ 
has the same cohomology class. 
\qed 

To understand the cocycle $\zeta$ associated to an affine connection, 
we first describe it in a local chart. 

\Remark II.2. If $U$ is an open subset 
of the vector space $V$, then any affine connection 
on $U$ is given by 
$$ \nabla_X Y = dY\cdot X + \Gamma(X,Y)  $$
with a $(2,1)$-tensor $\Gamma$. We then have 
$$ \eqalign{ 
{\cal L}_X (\nabla - \Gamma)(Y,Z) 
&= [X,dZ(Y)] - dZ([X,Y]) - d[X,Z](Y) \cr 
&= d^2 Z(X,Y) + dZ(dY(X)) - dX(dZ(Y)) \cr
&\ \ \ \ - dZ(dY(X)) + d Z(dX(Y)) - d(dZ(X) - dX(Z))(Y) \cr 
&= d^2 Z(X,Y) - dX(dZ(Y)) + d Z(dX(Y)) \cr
&\ \ \ \ 
- d^2 Z(Y,X) - dZ(dX(Y)) + d^2X(Y,Z) + dX(dZ(Y)) \cr 
&= d^2X(Y,Z). \cr} $$ 
This means that 
$$ {\cal L}_X \nabla = d^2 X + {\cal L}_X \Gamma. $$

Since the Lie derivative of any symmetric tensor is symmetric, we see that 
${\cal L}_X\nabla$ is symmetric if and only if ${\cal L}_X \Gamma$ is symmetric, 
which is the case if $\Gamma$ is symmetric, and this in turn means that 
$\nabla$ is torsion free. 
\qed

As a consequence of the preceding discussion, we can associate to any smooth manifold 
a canonical cohomology class $[\zeta] \in H^1_c({\cal V}_M, 
\Omega^1(M,\End(TM)))$ 
(cf.\ [Ko74]). That this class is always non-zero can be seen in local coordinates 
by showing that there is no $(1,2)$-tensor $\Gamma$ with $d^2 X = {\cal L}_X\Gamma$ 
for all smooth vector fields on a $0$-neighborhood in $\R^N$. 
For constant vector fields $X$, the preceding relation means that $\Gamma$ is constant 
and for linear vector fields $X$ we see that $\Gamma \: \R^N \otimes \R^N \to \R^N$ 
should be $\GL_N(\R)$-equivariant, so that 
$\lambda^2 \Gamma(v,w) =\Gamma(\lambda v,\lambda w) = \lambda\Gamma(v,w)$ 
for each $\lambda \in \R$ yields $\Gamma = 0$. 

\Example II.3. Let $V := \R^N$ and assume that $M$ is parallelizable of 
dimension~$N$. 
Then there exists some 
$\kappa \in \Omega^1(M,V)$ such that each $\kappa_m$ is invertible, i.e., 
$\kappa$  defines a trivialization 
of the tangent bundle of $M$ via $M \times V \to TM, 
(m,v) \mapsto \kappa_m^{-1}v$. 
Then, for each $X \in {\cal V}_M$, ${\cal L}_X\kappa$ also in an element of 
$\Omega^1(M,V)$, and since all the linear maps $\kappa_m$ are invertible, 
it can be written as 
${\cal L}_X.\kappa = -\theta(X)\cdot \kappa$ for some smooth function 
$\theta(X) \in C^\infty(M,\gl(V))$. 
Clearly, $C^\infty(M,\gl(V)) \cong C^\infty(M,\R) \otimes \gl(V)$ 
is a Lie algebra with respect to the pointwise bracket and ${\cal V}_M$ acts 
naturally by derivations. 

We claim that $\theta$ is a crossed homomorphism: 
$$ \eqalign{ 
\theta([X,Y]) \cdot \kappa 
&= - {\cal L}_{[X,Y]} \kappa = - {\cal L}_X({\cal L}_Y \kappa) + {\cal L}_Y({\cal L}_X\kappa) \cr 
&= {\cal L}_X(\theta(Y)\cdot \kappa) - {\cal L}_Y(\theta(X)\cdot\kappa) \cr 
&= ({\cal L}_X\theta(Y) - {\cal L}_Y\theta(X))\kappa 
+ \theta(Y){\cal L}_X\kappa  - \theta(X){\cal L}_Y\kappa \cr 
&= ({\cal L}_X\theta(Y) - {\cal L}_Y\theta(X))\kappa 
+ (\theta(Y)\theta(X) -\theta(X)\theta(Y))\kappa \cr 
&= ({\cal L}_X\theta(Y) - {\cal L}_Y\theta(X) + [\theta(X),\theta(Y)])\kappa.\cr} $$

\msk 

Now let 
$$ \nabla_X Y := \kappa^{-1}(X.\kappa(Y)) $$
denote the corresponding affine connection. Then 
$$ \kappa(\nabla_X Y) = X.\kappa(Y) = - \theta(X)\kappa(Y) + \kappa([X,Y]) \leqno(2.1) $$
is an immediate consequence of the definition of $\theta$. 
Moreover, the map 
$$ \tilde\kappa \: \Gamma(\End(TM)) \to C^\infty(M,\gl(V)), \quad \tilde\kappa(\phi)(m) = 
\kappa_m \circ \phi_m \circ \kappa_m^{-1} $$
satisfies 
$$ \tilde\kappa \circ {\cal L}_X\nabla = - d\theta(X) \in \Omega^1(M,\gl(V)).\leqno(2.2) $$
In fact, we have 
$$ \eqalign{ 
\kappa(({\cal L}_X\nabla)(Y)(Z)) 
&= \kappa([X,\nabla_Y Z] - \nabla_{[X,Y]}Z - \nabla_Y[X,Z]) \cr
&= \kappa(\nabla_X \nabla_Y Z) + \theta(X)\kappa(\nabla_Y Z) - [X,Y].\kappa(Z) - Y.\kappa([X,Z]) \cr 
&= XY.\kappa(Z) + \theta(X)(Y.\kappa(Z)) - [X,Y].\kappa(Z) - Y.(X.\kappa(Z)+\theta(X)\cdot \kappa(Z)) \cr 
&= \theta(X)(Y.\kappa(Z)) - Y.(\theta(X).\kappa(Z)) \cr 
&=  -Y.(\theta(X))\cdot \kappa(Z) = - d\theta(X)(Y)\cdot \kappa(Z), \cr } $$
and this calculation shows that 
$\tilde\kappa(({\cal L}_X\nabla)(Y)) = - d\theta(X)(Y)$
for each vector field $Y$, which is (2.2). 
\qed

\Remark II.4. For any affine connection $\nabla$, the operator 
$\eta(X)  := \nabla_X - \ad X$ 
on ${\cal V}_M$ is ${\cal F}_M$-linear, hence defines 
a section of $\End(TM)$. We thus obtain a map 
$$ \eta \: {\cal V}_M \to \Gamma(\End(TM)), \quad X \mapsto \nabla_X - \ad X. $$
The space $\Gamma(\End(TM))$ carries a natural associative algebra structure, 
hence in particular the structure of a 
Lie algebra (the gauge Lie algebra of the vector bundle $TM$) 
and ${\cal V}_M$ acts via the Lie derivative by derivations. 
The map $\eta$ is a crossed homomorphism for this structure if and only if 
$X \mapsto \nabla_X$ defines a representation of ${\cal V}_M$ on itself, i.e., 
if $\nabla$ is a flat connection (Remark~I.4). 
This is in particular the case if $\nabla$ is obtained 
from a trivialization $\kappa \in \Omega^1(M,V)$ 
of the tangent bundle. In the latter case, (2.1) 
implies that 
$$\tilde\kappa(\eta(X)) = - \theta(X). $$
In view of $\tilde\kappa(\nabla_X) = {\cal L}_X$, we now see that 
$$ \tilde\kappa(\ad X) = \tilde\kappa({\cal L}_X) = {\cal L}_X + \theta(X), $$
which is a representation of ${\cal V}_M$ on $C^\infty(M,V)$ if $\theta$ 
is a crossed homomorphism (Remark~I.4). 

We likewise see for the action of ${\cal V}_M$ on $\Gamma(\End(TM))$ and 
$C^\infty(M,\gl(V))$ that 
$\tilde\kappa(\phi) := \kappa \circ \phi \circ \kappa^{-1}$
leads to 
$$ \tilde\kappa \circ {\cal L}_X = ({\cal L}_X + \ad(\theta(X)) \circ \tilde\kappa, $$
i.e., the representation of ${\cal V}_M$ on $\Gamma(\End(TM))$ by the Lie derivative 
is transformed by $\tilde\kappa$ into the representation given by 
${\cal L}_X + \ad(\theta(X))$ on $C^\infty(M,\gl(V))$. 
Similarly, the map 
$$ \tilde\kappa_1 \: \Omega^1(M,\End(TM)) \to \Omega^1(M,\gl(V)), \quad 
\omega \mapsto \tilde\kappa \circ \omega $$
intertwines the Lie derivative on the left hand side with the 
representation on the right hand side given by ${\cal L}_X + \ad(\theta(X))$. 
\qed

\Definition II.5. Next 
we use the cocycle $\zeta$, associated to the affine connection 
$\nabla$, to define $k$-cocycles $\Psi_k \in Z^k_c({\cal V}_M, \Omega^k_M)$ 
depending on second order partial derivatives of vector fields. 
Let $V := \R^N$. For each $p \in \N$, we have a polynomial of 
degree $p$ on $\gl(V)$, invariant under conjugation, 
given by $A \mapsto \Tr(A^p)$. The corresponding invariant symmetric 
$p$-linear map is given by 
$$ \beta(A_1,\ldots, A_p) = \sum_{\sigma \in S_p} \Tr(A_{\sigma(1)}\cdots 
A_{\sigma(p)}), $$
and we consider it as a linear $\GL(V)$-equivariant map 
$\End(V)^{\otimes p} \to \R,$
where $\GL(V)$ acts trivially on $\R$. This $\GL(V)$-equivariant map leads to a 
vector bundle map 
$$ \beta_M \: \End(TM)^{\otimes p} \to M \times \R, $$
where $M \times \R$ stands for the trivial vector bundle with fiber $\R$. 
On the level of bundle-valued differential forms, 
this in turn yields an alternating $p$-linear map 
$$ \beta_M^1 \: \Omega^1(M,\End(TM))^{p} \to \Omega^p(M,\R) = \Omega^p_M. $$
To see that this map is ${\cal V}_M$-equivariant, we note that 
in local coordinates we have on an open subset $U \subeq V\cong 
\R^N$ the corresponding linear map 
$$ \Omega^1(U,\End(TU))^{\otimes p} 
\cong \Omega^1(U,\R)^{\otimes p} \otimes \End(V)^{\otimes p} \to \Omega^p(U,\R), $$
acting on $\Omega^1(U,\R)^{\otimes p}$ as the $p$-fold exterior product and 
on $\End(V)^{\otimes p}$ as $\beta$. 
If $\phi \: U_1 \to U_2$ is a local diffeomorphism, and 
$\alpha \in \Omega^1(U_2, \End(T(U_2)))$, then 
$\phi^*\alpha \in \Omega^1(U_1, \End(T(U_1)))$ is given by 
$$ (\phi^*\alpha)_m(v) := T_m(\phi)^{-1} \circ 
\alpha_{\phi(m)}(T_m(\phi)v) \circ T_m(\phi). $$ 
From this it follows that 
$\beta_M^1$ is equivariant under diffeomorphisms, and hence with respect to the 
infinitesimal action of vector fields. 

We conclude that we can use $\beta^1_M$ to multiply Lie algebra cocycles 
(cf.\ [Fu86], App.~F in [Ne04]). In 
particular, this leads with the cocycle $\zeta$ from Lemma~II.1 
for each $k \in \N$ to a Lie algebra cocycle 
$$ \Psi_k \in Z^k_c({\cal V}_M, \Omega^k_M), $$
defined by 
$$ \Psi_k (X_1, \ldots, X_k) 
:= (-1)^k \beta_M^1(\zeta(X_{\sigma(1)}), \ldots, \zeta(X_{\sigma(k)})) $$
(cf.\ [Fu86]). 

For any other affine connection $\nabla'$, the corresponding 
cocycle $\zeta'$, and the associated cocycles $\Psi_k'$, we first note that 
the difference $\zeta' - \zeta$ is a coboundary,  and since 
products of cocycles and coboundaries are coboundaries, 
$\Psi_k' - \Psi_k$ is a coboundary. 
Hence its cohomology class in $H^k_c({\cal V}_M, \Omega^k_M)$ 
does not depend on the choice of the affine connection. 
\qed

\Remark II.6. If the connection $\nabla$ is defined by a trivialization of $TM$ 
as in Example~II.3, then $\zeta(X)$ is transformed into $-d\theta(X) \in 
\Omega^1(M,\gl(V))$, which leads to 
$$ \Psi_k (X_1, \ldots, X_k) = \sum_{\sigma\in S_k} \sgn(\sigma)
\Tr \left( d\theta(X_{\sigma(1)}) \wedge \ldots \wedge d\theta(X_{\sigma(k)}) \right). 
\qeddis 

Next we construct several equivariant cocycles on a gauge Lie algebra. 

\Theorem II.7. Let $\g$ be a finite-dimensional Lie algebra,
$\F_M \otimes \g$ be the corresponding gauge Lie algebra, 
consider $\Omega^p_M$ as a trivial module of $\F_M \otimes \g$, 
and let $\rho \: \g \to \gl(V)$ 
be a finite-dimensional representation. 
\litem{(a)} By $\F_M$-linear extension, we get a morphism of cochain complexes, 
$C^\bullet(\g,\R) \rightarrow \linebreak C^\bullet_\equ (\F_M \otimes \g, \F_M)$, where equivariance 
refers to the action of ${\cal V}_M$, 
which gives a homomorphism of cohomology algebras 
$$ H^\bullet (\g,\R) \rightarrow H^\bullet_\equ (\F_M \otimes \g, \F_M) .$$

\litem{(b)} The following expression is a ${\cal V}_M$-equivariant cocycle 
with values in $k$-forms,
$\psi_k \in \linebreak Z^k_\equ (\F_M \otimes \g, \Omega^k_M)$:
$$ \psi_k (f_1 \otimes x_1, \ldots, f_k \otimes x_k)
= \Big(\sum_{\sigma\in S_k} \Tr \left( \rho(x_{\sigma(1)}) \cdots \rho(x_{\sigma(k)}) \right)\Big)
df_1 \wedge \cdots \wedge df_k ,$$
where $f_1, \ldots, f_k \in \F_M$, $x_1, \ldots, x_k \in\g$. 

\litem{(c)} The following expression is an equivariant cocycle with values in 
$\oline\Omega^{k-1}_M$, $k \geq 1$:
$$\bp_k (f_1 \otimes x_1, \ldots, f_k \otimes x_k)
= \Big(\sum_{\sigma\in S_k} \Tr \left( \rho(x_{\sigma(1)}) \cdots \rho(x_{\sigma(k)}) \right)\Big)
[f_1 df_2 \wedge \ldots \wedge df_k],$$
where $f_1, \ldots, f_k \in \F_M$, $x_1, \ldots, x_k \in\g$. 

\Proof. Verification of part (a) is straightforward. 
Assertions (b) and (c) easily follow from Proposition~A.3, which is first used  with 
$\Omega^p = \Omega^p_M$ to obtain the cocycle property on 
${\cal F}_M \otimes \gl(V)$, and then we pull it back by the homomorphism 
$\id_{{\cal F}_M} \otimes \rho \: {\cal F}_M \otimes \g \to {\cal F}_M\otimes \gl(V)$. 
\qed

\Definition II.7. Combining Theorems I.7 and II.7 with Renark~II.4, 
we get for each trivialization of the tangent bundle $T(M)$ 
of an $N$-dimensional manifold $M$ 
by pulling back with the corresponding crossed homomorphism 
$\theta \: {\cal V}_M \to \F_M \otimes \gl(V)$, $V = \R^N$, 
the following cocycles 
$$ \Phi_k \in Z_c^{2k-1} ({\cal V}_M, \F_M), \quad 
\Phi_k (X_1, \ldots, X_{2k-1}) = \sum_{\sigma\in S_{2k-1}} \sgn(\sigma)
\Tr \left( \theta(X_{\sigma(1)}) \cdots \theta(X_{\sigma(2k-1)}) \right)$$
and $\Bp_k \in Z^k_c ({\cal V}_M, \oline\Omega^{k-1}_M)$, defined by 
$$ \Bp_k (X_1, \ldots, X_k) = \sum_{\sigma\in S_k} \sgn(\sigma)
[\Tr \left( \theta(X_{\sigma(1)}) d\theta(X_{\sigma(2)}) \wedge \ldots \wedge d\theta(X_{\sigma(k)}) \right)].$$ 
Note that we have for each $k\geq 1$ the relation 
$$ d \circ \oline\Psi_k = \Psi_k. \leqno(2.3)$$ 
\qed

\sectionheadline{III. Cohomology 
of smooth vector fields with values in differential forms}

In this section, we first recall Tsujishita's Theorem ([Tsu81], Thm.~5.1.6, see also 3.3.4),  
describing the continuous cohomology of ${\cal V}_M$ with the values in differential forms. 

To explain its statement, we recall that each compact manifold 
$M$ satisfies Tsujishita's condition (F): $H^\bullet_M$ is finite-dimensional 
and, for each $p \in \N$, the subspace 
$Z^p_M$ of closed $p$-forms in $\Omega^p_M$ has a closed complement.  
For compact oriented manifolds, the 
latter assertion is a direct consequence of the Hodge Decomposition Theorem 
(cf.\ [AMR83, Thm.~7.5.3]) and the general case can be reduced to this one via 
the orientable $2$-sheeted covering. Moreover, 
Beggs shows in [Beg87] (Theorems~4.4 and 6.6) that for all finite-dimensional 
paracompact smooth manifolds the space 
$Z^p_M$ of closed $p$-forms in $\Omega^p_M$ has a closed complement, so that 
condition (F) is equivalent to 
all cohomology spaces $H^p_M$, resp., the cohomology algebra 
$H^\bullet_M$, being finite-dimensional.  

The following theorem is a key result of this paper 
because it provides a bridge between 
Tsujishita's results an the cocycles $\Psi_k$. It will be proved in Appendix C. 
 
\Theorem III.1. Let $M$ be an $N$-dimensional manifold 
for which $H^\bullet_M$ is finite-dimensional. Then  
$H^q_c({\cal V}_M, \Omega^p_M)$ vanishes if $p > q$ and 
for $p,n \in \N_0$,  we have 
$$H^{p+n}_c({\cal V}_M, \Omega^p_M) \cong 
H^n_c ({\cal V}_M, \F_M) \otimes E_p, $$
where 
$$E_p := \span \Big\{ [\Psi_1]^{m_1}[\Psi_2]^{m_2} \cdots [\Psi_p]^{m_p} \: 
\sum_{j=1}^p j m_j = p \Big\} \subeq H^p_c({\cal V}_M, \Omega^p_M). $$
In particular, 
$H^\bullet_c({\cal V}_M, \Omega^\bullet_M)$ is a free 
$H^\bullet ({\cal V}_M, \F_M)$-module, generated by the non-zero 
products of the classes $[\Psi_k]$, $k = 1,\ldots, N$. 
\qed

From the preceding theorem, we immediately derive with 
$H^0({\cal V}_M, {\cal F}_M) = H^0_M \cong \R$: 

\Corollary III.2. 
\litem{(a)} If $p < k$, then $H^p ({\cal V}_M, \Omega_M^k) = 0$. 
\litem{(b)} If $M$ is connected, then 
$H^p ({\cal V}_M, \Omega^p_M) \cong E_p$. 
\qed

\Lemma III.3.  Let 
$M$ be a compact connected orientable manifold and 
$\mu$ a volume form on $M$. Then 
$$ \div \: {\cal V}_M \to {\cal F}_M, \quad {\cal L}_X\mu = \div(X) \mu $$
defines a continuous Lie algebra cocycle and furthermore, each closed $1$-form 
$\alpha \in Z^1_M$ defines a Lie algebra cocycle 
${\cal V}_M \to {\cal F}_M, X \mapsto \alpha(X).$ 
These cocycles span the cohomology space 
$$H^1_c({\cal V}_M, {\cal F}_M) \cong H^1_M \oplus \R[\div]. $$
If, in addition, the volume form is defined by a trivialization of 
the tangent bundle $T(M)$, then $-\div = \Phi_1 = \oline\Psi_1$. 

\Proof. In view of [Fu86, Th.~II.4.11], we only have to verify the last part 
(cf.\ also [FL80] for the local cohomology of ${\cal V}_M$ with values 
in $\Omega^p_M$). 
So let us assume that the volume form can be written as 
$$ \mu = \kappa_1 \wedge \ldots \wedge \kappa_N, $$
where the $\kappa_i$ are the components of a trivializing 
$1$-form $\kappa \in \Omega^1(M,\R^N)$. 
From ${\cal L}_X\kappa = - \theta(X)\kappa$ we derive 
${\cal L}_X\kappa_j = - \sum_{i = 1}^N \theta(X)_{ji} \kappa_i$ 
and further 
$$ \eqalign{ {\cal L}_X\mu 
&= \sum_{i=1}^N \kappa_1 \wedge \ldots \wedge 
{\cal L}_X \kappa_i \wedge \ldots \wedge \kappa_N 
= -\Big(\sum_{i=1}^N \theta(X)_{ii}\Big) \kappa_1 \wedge \ldots \wedge \kappa_N 
= -\Tr(\theta(X)) \cdot \mu. \cr}$$
This proves that $\div X = - \Phi_1(X) = - \oline \Psi_1(X).$ 
\qed

\Remark III.4. We take a closer look at the second cohomology of the ${\cal V}_M$-modules 
$\Omega^p_M$. 

From Corollary~III.2(a) we know that $H^2_c({\cal V}_M, \Omega_M^p)$ 
vanishes for $p > 2$, and for $p = 2$ the space 
$$ H^2_c({\cal V}_M, \Omega^2_M) = \span \{[\Psi_1^2], [\Psi_2]\} $$
is $2$-dimensional. For $p = 1$ we also know that 
$$ H^2_c({\cal V}_M, \Omega^1_M) 
= H^1_c({\cal V}_M, {\cal F}_M).[\Psi_1] 
= H^1_M\cdot [\Psi_1]\oplus \R [\oline\Psi_1 \wedge \Psi_1] $$ 
is of dimension $b_1(M) + 1$ (Lemma III.3). 
The most intricate case is $p =0$, i.e., $H^2_c({\cal V}_M, {\cal F}_M)$ 
(cf.\ Theorem~IV.8 below. 
\qed

\Example III.5. Let $M = G$ be an abelian Lie group of dimension $N$ 
and $\kappa_1, \ldots, \kappa_N$ a basis of the space of invariant 
$\R$-valued $1$-forms on $G$. Then 
$\kappa := (\kappa_1,\ldots, \kappa_N)$ is a closed trivializing 
$1$-form (the Maurer--Cartan form of $G$). 
Let $X_1,\ldots, X_N \in {\cal V}(G)$ be the dual basis of left invariant 
vector fields. Writing $X \in {\cal V}(G)$ as 
$X = \sum_{i=1}^N f_i X_i$, we obtain 
$$ {\cal L}_X \kappa_i = d(i_X \kappa_i) = df_i. $$
Therefore $\theta(X) \in C^\infty(G,\gl_N(\R))$ is given by 
$\theta(X)_{ij} = -df_i(X_j)$, i.e., $-\theta(X)$ is the Jacobian 
matrix of $X$ with respect to the basis $X_1,\ldots, X_N$. 
\qed

\sectionheadline{IV. Cohomology with values in differential forms modulo exact forms} 

Let $M$  be a parallelizable connected compact manifold of dimension $N$. 
Since the subspace 
$H^k_M$  of $\oline\Omega^k_M$ is a trivial ${\cal V}_M$-module, we derive from 
Corollary~D.5 that 
$$ H_c^k({\cal V}_M, H^m_M) = \0 \quad \hbox{ for } \quad 0 < k \leq N, m \in \N_0. \leqno(4.1) $$
Therefore the short exact sequence of ${\cal V}_M$-modules 
$$ \0 \to H^m_M \to \oline\Omega^m_M \sssmapright{d} B^{m+1}_M \to \0 $$
induces a long exact sequence in cohomology, hence leads to isomorphisms 
$$ d_* \: H_c^k({\cal V}_M, \oline\Omega^m_M) \to H_c^k({\cal V}_M, B^{m+1}_M) \quad \hbox{ for  } \quad 
m \in \N_0, 0 < k \leq N-1 \leqno(4.2a) $$
and an exact sequence 
$$ \0 \to H_c^N({\cal V}_M, \oline\Omega^m_M) \ssmapright{d_*} 
H_c^N({\cal V}_M, B^{m+1}_M) \to 
 H_c^{N+1}({\cal V}_M, H^m_M) 
\cong  H_c^{N+1}({\cal V}_M, \R) \otimes H^m_M.\leqno(4.2b) $$ 
In particular, we get 
$$ \dim H_c^N({\cal V}_M, \oline\Omega^m_M) \leq \dim H_c^N({\cal V}_M, B^{m+1}_M). $$

\Lemma IV.1. {\rm([Ne06a, Lemma~23], [HS53])} 
For any smoothly paracompact manifold $M$ we have 
$$ H_c^0({\cal V}_M, \oline\Omega^p_M) = (\oline\Omega^p_M)^{{\cal V}_M} = H^p_M. 
\qeddis 

\Lemma IV.2. If $M$ is parallelizable, then for 
each $k \in \N$ the exterior differential $d$ induces a surjective map 
$$ d_* \: H^k_c({\cal V}_M, \oline\Omega^{k-1}_M) \to 
H^k_c({\cal V}_M, \Omega^k_M), \quad [\alpha] \mapsto [d \circ \alpha], $$
and the natural map 
$$ H^k_c({\cal V}_M, B^k_M) \to H^k_c({\cal V}_M, \Omega^k_M) $$
is surjective. 

\Proof. First we recall from Definition~II.7 that $d\circ \oline\Psi_k = \Psi_k.$ 

Next we observe that the exterior product of differential forms induces maps 
$$ Z^p_M \times B^q_M \to B^{p+q}_M, 
\quad (\alpha,\beta) \mapsto \alpha \wedge \beta, $$
and hence maps 
$$ Z^p_M \times \oline\Omega^q_M \to \oline\Omega^{p+q}_M, 
\quad (\alpha,[\beta]) \mapsto [\alpha \wedge \beta]. $$
Since $\Psi_j$ has values in $B^1_M \wedge \ldots \wedge B^1_M \subeq B^j_M$, 
each product 
$\Psi_{j_1} \wedge \ldots \wedge \Psi_{j_r}$
has values in $B^{j_1+\ldots + j_r}_M$, so that 
$$ \oline\Psi_{j_1} \wedge \Psi_{j_2} \wedge\ldots \wedge \Psi_{j_r} $$
is a well-defined cocycle with values in 
$\oline\Omega_M^{j_1+\ldots + j_r}$, satisfying 
$$ d \circ \Big(\oline\Psi_{j_1} \wedge \Psi_{j_2} \wedge\ldots \wedge \Psi_{j_r}\Big)  
=  (d \circ \oline\Psi_{j_1}) \wedge \Psi_{j_2} \wedge\ldots \wedge \Psi_{j_r}  
=  \Psi_{j_1} \wedge \Psi_{j_2} \wedge\ldots \wedge \Psi_{j_r}. $$
For $j_1 + \ldots + j_r = k$, this implies that the image of $d$ contains 
all products 
$[\Psi_{j_1} \wedge \Psi_{j_2} \wedge\ldots \wedge \Psi_{j_r}]$ 
with $\sum j_i = k$, and, in view of Theorem~III.1, these products 
span $H^k_c({\cal V}_M, \Omega_M^k)$. 

This proves the first part, and the second part of the assertion is an immediate 
consequence of the first one. 
\qed

In view of 
$$ H_c^q({\cal V}_M, \Omega^p_M) = \0 \quad \hbox{ for } \quad q < p\leqno(4.3)  $$ (Corollary~III.2(a)), the long exact cohomology sequence associated to the short exact sequence  
$$ \0 \to B^p_M \to \Omega^p_M \to \oline\Omega^p_M \to \0 $$
(which splits topologically by condition (F)), 
leads to isomorphisms 
$$ H_c^q({\cal V}_M, \oline\Omega^p_M) \cong H_c^{q+1}({\cal V}_M, B^p_M) \quad \hbox{ for } \quad q < p -1. 
\leqno(4.4a) $$
For $q = p-1$, Lemma~IV.2 leads to an exact sequence 
$$ \0= H_c^{p-1}({\cal V}_M, \Omega^p_M) \to H_c^{p-1}({\cal V}_M, \oline\Omega^{p}_M) 
\to H_c^p({\cal V}_M, B^p_M) \onto H_c^p({\cal V}_M, \Omega^p_M) \to \0.\leqno(4.4b) $$

From [Ne06a, Prop.~6], we recall: 
\Lemma IV.3. For each closed $(p+q)$-form $\omega \in \Omega^{p+q}_M$,  
the prescription 
$$ \omega^{[p]}(X_1,\ldots, X_p) := [i_{X_p} \ldots i_{X_1} \omega] \in \oline\Omega^q_M $$
defines a continuous $p$-cocycle in $Z^p_c({\cal V}_M,\oline\Omega^q_M)$. 
\qed

\subheadline{Calculating $H_c^2({\cal V}_M, \oline\Omega^m_M)$} 

In this subsection we address the problem to determine the second cohomology 
space $H_c^2({\cal V}_M, \oline\Omega^m_M)$ for a parallelizable manifold $M$  
of dimension~$N$. 

We start by observing that we have the following isomorphisms 
$$ \leqalignno{ 
H_c^2({\cal V}_M, \oline\Omega^{m}) 
&\cong H_c^2({\cal V}_M, B_M^{m+1}) \quad \hbox{ for } N \geq 3 \quad \hbox{ by } \quad (4.2a)\cr 
&\cong H_c^1({\cal V}_M, \oline\Omega^{m+1}_M) \quad \hbox{ for } m > 1 \quad \hbox{ by } \quad (4.4a)&(4.5)\cr 
&\cong H_c^1({\cal V}_M, B_M^{m+2})  \quad \hbox{ for } N \geq 2 \quad \hbox{ by } \quad (4.2a)\cr 
&\cong H_c^0({\cal V}_M, \oline\Omega^{m+2}_M)\quad \hbox{ for each } m \quad \hbox{ by } \quad 
(4.4a) \cr  
&= H^{m+2}_M \quad \hbox{ by Lemma~IV.1.}\cr  } $$

We thus obtain 
\Proposition IV.4. {\rm(a)} For $N \geq 3$ and $m \geq 2$ the map 
$$ H^{m+2}_M \to H_c^2({\cal V}_M, \oline\Omega^m_M), \quad 
[\omega] \mapsto [\omega^{[2]}] $$
is an isomorphism. 

{\rm(b)} For $N \geq 2$ and $m \geq 1$ the map 
$$ H^{m+1}_M \to H_c^1({\cal V}_M, \oline\Omega^m_M), \quad 
[\omega] \mapsto [\omega^{[1]}] $$
is an isomorphism. 

\Proof. In view of (4.5), we only have to see that the isomorphisms are 
implemented by the maps $\omega \mapsto \omega^{[2]}$, resp., $\omega \mapsto \omega^{[1]}$. 

(a) We start with a closed $(m+2)$-form $\omega$, representing an element 
of $H^{m+2}_M$, which we consider a ${\cal V}_M$-fixed element in 
$\oline\Omega^{m+2}_M$. A corresponding $1$-cocycle with values in 
$B^{m+2}_M$ is given by $\alpha(X) := {\cal L}_X\omega = d(i_X \omega) = 
(d \circ \omega^{[1]})(X)$. This already shows that the corresponding 
$\oline\Omega^{m+1}_M$-valued $1$-cocycle is $\omega^{[1]}$ (cf.\ Lemma~IV.3). 

To distinguish it from the exterior differential $d$, we write 
$d_{{\cal V}(M)}$ for the Chevalley--Eilenberg differential on 
$C^\bullet({\cal V}_M, \Omega^\bullet_M)$. 
To see the corresponding $2$-cocycle with values in $B_M^{m+1}$, we recall 
from [Ne06a, Lemma~5] that for $\tilde\omega^{[1]}(X) := i_X \omega$ we have 
$$ (d_{{\cal V}_M}\tilde\omega^{[1]})(X,Y) = d(i_X i_Y \omega) 
= -(d \circ \omega^{[2]})(Y,X). \leqno(4.6) $$
Hence $-\omega^{[2]} \in Z^2_c({\cal V}_M, \oline\Omega^m_M)$ is the 
$2$-cocycle corresponding to $[\omega]$ under the isomorphisms 
$H^{m+2}_M \cong H^2_c({\cal V}_M,\oline\Omega^m)$ in (4.5). 

(b) is proved as in (a). 
\qed

Now we turn to the cases not covered by Proposition~IV.4. 
The following proposition was the original motivation for the present paper: 

\Theorem IV.5. For $N \geq 2$ we have 
$$ H^2_c({\cal V}_M,\oline\Omega^1_M) \cong H^3_M \oplus \R [\oline\Psi_1 \wedge \Psi_1] 
\oplus \R[\oline\Psi_2], $$
where $H^3_M$ embeds via the map $[\omega] \mapsto [\omega^{[2]}]$. 

\Proof. First we assume that $N \geq 3$. 
For $m = 1$, we first get from (4.2a) that 
$H^2_c({\cal V}_M,\oline\Omega^1_M) \cong H_c^2({\cal V}_M,B_M^2)$
and further from (4.4b) the exact sequence 
$$ \0 \to H_c^1({\cal V}_M, \oline\Omega_M^2)\into H_c^2({\cal V}_M, B_M^2) 
\to H_c^2({\cal V}_M,\Omega^2_M) = \R [\Psi_1^2] \oplus 
\R [\Psi_2], $$
where the inclusion on the left maps $[\omega^{[1]}]$ to $[\omega^{[2]}]$ 
(cf.\ the proof of Proposition~IV.4). 
Since we know from Proposition~IV.4 that the map 
$$H^3_M \to H_c^1({\cal V}_M,\oline\Omega^2_M), \quad [\omega] \mapsto 
[\omega^{[1]}] $$ 
is an isomorphism, the assertion follows 
 because the classes $[\Psi_1^2]$ and $[\Psi_2]$ in $H_c^2({\cal V}_M, \Omega^2_M)$ 
are 
contained in the image of $H_c^2({\cal V}_M,B_M^2)$ (Lemma~IV.2) and 
satisfy 
$$ [d \circ (\oline\Psi_1 \wedge \Psi_1)] = [\Psi_1^2]
\quad \hbox{ and } \quad 
[d \circ \oline\Psi_2] = [\Psi_2]. $$

Now we consider the case $N = 2$. 
For $m = 1$ we get from (4.2a) and Proposition IV.4 
$$ H_c^1({\cal V}_M,\oline\Omega_M^1) \cong H_c^1({\cal V}_M,B_M^2) \cong H^2_M 
\quad \hbox{ and } \quad 
 H_c^1({\cal V}_M,\oline\Omega_M^2) \cong 
 H_c^1({\cal V}_M,B^3_M) = \{0\}. $$
Therefore the exact sequence 
$B_M^2 \into \Omega_M^2 \to \oline\Omega_M^2$ leads to an exact sequence 
$$ \0= H_c^1({\cal V}_M, \oline\Omega_M^2) \to H_c^2({\cal V}_M,B_M^2) \to H_c^2({\cal V}_M, \Omega_M^2) \to H_c^2({\cal V}_M, \oline\Omega_M^2)=\0, $$
showing that the map 
$$H_c^2({\cal V}_M,B_M^2) 
\to H_c^2({\cal V}_M,\Omega_M^2) \cong \R[\Psi_1^2] \oplus \R[\Psi_2]$$ 
is an isomorphism. 

On the other hand, (4.2b) yields an embedding 
$$ d\circ \: H_c^2({\cal V}_M,\oline\Omega_M^1) \into H_c^2({\cal V}_M, B_M^2) 
\cong H_c^2({\cal V}_M,\Omega_M^2) \cong \R[\Psi_1^2] \oplus \R[\Psi_2], $$ 
and Lemma IV.2 asserts that this embedding is surjective. 
This completes the proof. 
\qed

\Proposition IV.6. For $N = 2$ we have $H_c^2({\cal V}_M, \oline\Omega^m_M) = \0$ for $m \geq 2$ and 
$$ \dim H^2_c({\cal V}_M,\oline\Omega^1_M) = 2. $$

\Proof. For $N = 2$ we have $\Omega_M^m = 0$ for $m \geq 3$, so that $\oline\Omega_M^m$ vanishes 
in these cases. For $m = 2$ we get from (4.2b) an inclusion 
$H_c^2({\cal V}_M,\oline\Omega_M^2) \into H_c^2({\cal V}_M,B_M^3) = \0$. Therefore 
$H_c^2({\cal V}_M,\oline\Omega_M^m)$ vanishes for $m \geq 2$. The 
case $m = 1$ follows from Theorem~IV.5. 
\qed

So far we have covered all cases $m \geq 1$ for $N \geq 2$, and 
we are left with the case $m = 0$ or $N = 1$. 
First we recall the results for $N = 1$: 

\Example IV.7. (a) (cf.\ [FF01, Th.~30]) For $N =1$ and $M = \SS^1$ the algebra 
$H^\bullet({\cal V}_{\SS^1}, {\cal F}_{\SS^1})$ is a free graded commutative 
algebra generated by two elements $\alpha_1, \alpha_2$ of degree $1$ and one element 
$\beta$ of degree $2$. 
In particular 
$$ H^1_c({\cal V}_{\SS^1}, {\cal F}_{\SS^1}) = \R\alpha_1 \oplus \R \alpha_2 \quad 
\hbox{ and } \quad 
H^2_c({\cal V}_{\SS^1}, {\cal F}_{\SS^1}) = \R \alpha_1\alpha_2 \oplus \R \beta. $$
Furthermore $\oline\Omega^1_{\SS^1} = H^1_{\SS^1}$ is one-dimensional and a trivial 
${\cal V}_{\SS^1}$-module, so that 
$$ H^2({\cal V}_{\SS^1}, \oline\Omega^1_{\SS^1}) = 
H^2({\cal V}_{\SS^1}, \R) \cong \R $$
(Theorem~29 in [FF01], p.~195). 

(b) For $M = \R^N$ the cohomology algebra 
$H^\bullet({\cal V}_{\R^N}, {\cal F}_{\R^N})\cong H^\bullet(\gl_N(\R),\R)$ 
is the free exterior 
algebra generated by the classes $[\Phi_k]$, $k =1,\ldots, N$ 
(cf.\ Definition~II.7; [Tsu77], \S 4). 
\qed

Now we turn to the case $m = 0$. 

\Theorem IV.8. For $N \geq 2$ the map 
$$ H^2_M \oplus H^1_M \to H^2_c({\cal V}_M,{\cal F}_M), \quad 
([\alpha], [\beta]) \mapsto 
[\alpha + \beta \wedge \oline\Psi_1] $$
is a linear isomorphism. 

\Proof. In view of $\oline\Omega_M^0 = {\cal F}_M$, (4.2a/b) provide an 
embedding 
$$ d_* \: H_c^2({\cal V}_M, {\cal F}_M) \into H_c^2({\cal V}_M, B_M^{1}), \quad 
[\omega] \mapsto [d \circ \omega] $$
which is an isomorphism for $N > 2$. 
We also get from Proposition IV.4 for $N \geq 2$ an isomorphism:  
$$ H^2_M \to H_c^1({\cal V}_M, \oline\Omega_M^1), \quad [\omega] \mapsto 
[\omega^{[1]}]. $$
In view of Lemma~III.3, we have 
$$ H_c^2({\cal V}_M, \Omega_M^1) = H_c^1({\cal V}_M, {\cal F}_M) \cdot[\Psi_1] \cong H^1_M 
\oplus \R [\oline\Psi_1 \wedge \Psi_1]. $$ 

We now consider the exact sequence 
$$ \eqalign{ & H_c^1({\cal V}_M,B^1_M) 
\onto H_c^1({\cal V}_M,\Omega_M^1) 
\ssssmapright{0} 
H_c^1({\cal V}_M,\oline\Omega_M^1) \cong H^2_M \cr
&\ \ \ \ssmapright{-d_*} H_c^2({\cal V}_M,B_M^1) \sssmapright{f} 
H_c^2({\cal V}_M,\Omega_M^1) \sssmapright{g}
H_c^2({\cal V}_M,\oline\Omega_M^1), \cr} $$
in which we know all terms except the middle one, which we are interested in. 
The one-dimensional space $H_c^1({\cal V}_M,\Omega_M^1)$ is generated by 
$[\Psi_1],$ 
which has values in $B_M^1$, so that its image in $\oline\Omega_M^1$ vanishes. 
This yields the $0$-arrow in the upper row of the diagram (Lemma~IV.2). 
To calculate the connecting map 
$H_c^1({\cal V}_M,\oline\Omega_M^1) \cong H^2_M \to H_c^2({\cal V}_M,B_M^1)$, 
we recall from (4.6) that we have for each closed $2$-form $\omega$ the 
relation 
$$\delta([\omega^{[1]}]) = [d_{{\cal V}_M} \tilde\omega^{[1]}] 
= - [d \circ \omega^{[2]}] = - d_*[\omega]. $$ 

Let $\beta$ be a closed $1$-form on $M$, considered as a $1$-cocycle 
${\cal V}_M \to {\cal F}_M$. 
Then 
$$ \beta \wedge \oline\Psi_1 \in Z_c^2({\cal V}_M, {\cal F}_M) $$
satisfies 
$$ d \circ (\beta \wedge \oline\Psi_1) 
= (d\circ \beta) \wedge \oline\Psi_1 
+ \beta \wedge d\circ \oline\Psi_1 
= (d\circ \beta) \wedge \oline\Psi_1 
+ \beta \wedge \Psi_1. $$
The $1$-cochain 
$$\gamma \: {\cal V}_M \to \Omega^1_M, \quad 
\gamma(X) := \oline\Psi_1(X)\cdot \beta = (\oline\Psi_1 \cdot \beta)(X)$$
satisfies 
$$ d_{{\cal V(M)}}\gamma 
= (d_{{\cal V}(M)}\oline\Psi_1) \wedge \beta
- \oline\Psi_1 \wedge (d_{{\cal V}(M)}\beta)
= -\oline\Psi_1 \wedge d_{{\cal V}(M)}\beta. $$
Since the Cartan formula implies that 
$$  (d_{{\cal V}(M)}\beta)(X) = {\cal L}_X\beta 
= i_X(d\beta) + d(i_X\beta) = d(\beta(X)), $$
we obtain 
$d_{{\cal V(M)}}\gamma 
= -\oline\Psi_1 \wedge (d \circ \beta),$
showing that 
$(d\circ \beta) \wedge \oline\Psi_1 \in B_c^2({\cal V}_M,\Omega_M^1).$
This implies that 
$$ d_*[\beta \cdot \oline\Psi_1] = [d \circ (\beta \wedge \oline\Psi_1)] = [\beta \wedge \Psi_1]. $$

We also recall that 
$$ H_c^2({\cal V}_M,\Omega_M^1) = (H^1_M \wedge [\Psi_1]) \oplus 
\R (\oline\Psi_1 \wedge \Psi_1) \cong H^1_M \oplus \R. $$
From the preceding considerations, we see that the subspace 
$H^1_M \wedge [\Psi_1]$ lies in the range of $f$. 
Since $g([\oline \Psi_1 \wedge \Psi_1]) = [\Psi_1^2]$ is non-zero, 
we derive 
$$ \im(f) = \ker(g) =H^1_M \wedge [\Psi_1]. $$
From the exactness in $H^1_c({\cal V}_M,B^1_M)$ we now see that the 
natural map 
$$ H^2_M \oplus H^1_M \to H^2_c({\cal V}_M,B^1_M), \quad 
([\alpha],[\beta]) \mapsto 
[d \circ \alpha + \beta \wedge \Psi_1] 
= [d \circ (\alpha + \beta \cdot \oline\Psi)] $$
is a linear isomorphism. 
As the elements $\alpha + \beta \cdot \oline\Psi$ 
form ${\cal F}_M$-valued $2$-cocycles, we finally conclude 
that the inclusion $d_*$ is also surjective for  $N = 2$. 
This completes the proof. 
\qed

\Problem 1. Determine the spaces $H^2({\cal V}_M,\oline\Omega^1_M)$ for 
all connected smooth manifolds~$M$, without assuming that $M$ 
is parallelizable. 
\qed

\sectionheadline{Appendix A. Lie algebra cocycles with values in associative algebras} 

Let $A$ be a unital associative algebra 
and $A_L := (A,[\cdot,\cdot])$ be the underlying Lie 
algebra. We consider $A$ as an $A_L$-module with respect to the adjoint representation 
$x.y := [x,y]$. 

We write $m_A \: A \times A \to A$ for the product map and $b_A \: A \times A \to A$ for 
the commutator bracket and observe that both are $A_L$-invariant. 

We take a closer look at the Lie algebra complex $(C^\bullet(A_L, A),d)$, 
where we use the multiplication on $A$, which is $A_L$-equivariant, to define a 
multiplication on $C^\bullet(A_L,A)$ by 
$$ (\alpha \wedge \beta)(x_1,\ldots, x_{p+q}) 
:= {1\over p!q!} \sum_{\sigma \in S_{p+q}} 
\sgn(\sigma) 
m_A\big(\alpha(x_{\sigma(1)}, \ldots, x_{\sigma(p)}), 
\beta(x_{\sigma(p+1)}, \ldots, x_{\sigma(p+q)})\big) $$
for $\alpha \in C^p(A_L,A)$ and $\beta \in C^q(A_L,A)$ (cf.\ [Fu86, Sect.~I.3.2]). 
With 
$$ \Alt(\gamma)(x_1, \ldots, x_r) 
:= \sum_{\sigma \in S_r} \sgn(\sigma) \alpha(x_{\sigma(1)}, \ldots, x_{\sigma(r)}), $$
this means that 
$$ \alpha \wedge \beta 
= {1 \over p! q!} \Alt(\alpha \cdot \beta)
= {1 \over p! q!} \Alt(m_A \circ (\alpha,\beta)). $$
An easy induction implies that $n$-fold products in $C^\bullet(A_L, A)$ are given by 
$$ \alpha_1 \wedge \cdots \wedge \alpha_n 
= {1 \over \prod_{j = 1}^n p_j!} \Alt(\alpha_1 \cdots \alpha_n) \quad \hbox{ for } \quad 
\alpha_i \in C^{p_i}(A_L,A). $$
We thus obtain an associative differential graded algebra 
$(C^\bullet(A_L,A),d)$ (cf.\ [Ne04, App.~F]; see also [Tsu81, p.~30]). 
\msk 

Note that $\id_A \in C^1(A_L,A)$ is a $1$-cochain with 
$$ d_{A_L}(\id_A)(x,y) = x.y - y.x - [x,y] = [x,y] - [y,x] - [x,y] = [x,y], $$
so that $d_{A_L}(\id_A) = b_A$. It follows in particular that $b_A$ is a $2$-coboundary. 
We therefore obtain a sequences of coboundaries 
$$ b_A^n = d_{A_L}(\id_A \wedge b_A^{n-1}) \in B^{2n}(A_L, A) \subeq Z^{2n}(A_L,A). $$

\Lemma A.1. If $T \: A \to \z$ is a linear map vanishing on all commutators, then 
$T \circ b_A^k$ vanishes for each $k \in \N$, and if $\z$ is considered as a trivial 
$A_L$-module, then  
$$ \phi_k := T \circ (\id_A \wedge b_A^{k-1}) = T \circ (\id_A)^{2k-1} 
\in Z^{2k-1}(A_L,\z), $$
is a cocycle given by 
$$ \phi_k(x_1, \ldots, x_{2k-1}) = 
\sum_{\sigma \in S_{2k-1}} \sgn(\sigma) T(x_{\sigma(1)} \cdots x_{\sigma(2k-1)}). $$

\Proof. Arguing as in Remark~I.2 of [Ne06b], we see that for each $n \in \N_0$, we have 
$$ b_A^n= {1\over 2^n} \Alt(b_A \cdots b_A) 
= \Alt(m_A \cdots m_A) = \Alt(\id_A^{2n}) = (\id_A)^{2n}, $$
i.e., 
$$ b_A^n(x_1, \ldots, x_{2n}) 
= \sum_{\sigma \in S_{2n}} \sgn(\sigma) x_{\sigma(1)} \cdots x_{\sigma(2n)}. $$

Let $\tau$ be the cyclic permutation $\tau = (1\ 2\ \ldots\ 2n)$. 
Then $\sgn(\tau) = -1$, but, for each $\sigma \in S_{2n}$, we have  
$$ T(x_{\sigma\tau(1)} \cdots x_{\sigma\tau(2n)}) 
= T(x_{\sigma(2)} \cdots x_{\sigma(2n)}x_{\sigma(1)}) 
= T(x_{\sigma(1)}x_{\sigma(2)} \cdots x_{\sigma(2n)}). $$
This leads to $T \circ b_A^n = - T_A \circ b_A^n$, which implies that 
$T \circ b_A^n$ vanishes. 

Since $b_A$ is a cocycle, we now obtain 
$$ d\phi_k = d(T \circ (\id_A \wedge b_A^{k-1})) 
= T \circ (d(\id_A) \wedge b_A^{k-1})
= T \circ (b_A^k) = 0. 
\qeddis

The following theorem describes an important application of the preceding 
construction, namely that for $A = \gl_n(\R)$ we thus obtain a set of generators 
of the cohomology algebra ([Fu86, Th.~2.1.6]): 

\Theorem A.2. The cohomology algebra $H^\bullet(\gl_N(\R),\R)$ 
is generated by the cohomology classes of the cocycles 
$\phi_k \in Z^{2k-1}(\gl_N(\R),\R), k =1,\ldots, N,$
given by 
$$ \phi_k (x_1, \ldots, x_{2k-1} ) = \sum\limits_{\sigma\in S_{2k-1}} 
\sgn(\sigma) \Tr \left( x_{\sigma(1)} \ldots x_{\sigma(2k-1)} \right) \quad 
\hbox{ for } \quad x_1, \ldots, x_{2k-1} \in \gl_N(\R).
\qeddis

\subheadline{Lie algebra cocycles with values in differential forms} 

Let $(\Omega^\bullet,d)$ be a differential graded algebra 
whose multiplication is denoted $\alpha \wedge \beta$, 
$V$ a finite-dimensional vector space, and put 
$A := \Omega^0 \otimes \End(V)$ and $B := \Omega^\bullet \otimes \End(V)$. Then 
$B$ also carries the structure of a differential 
graded algebra with differential $d_B := d \otimes \id_{\End(V)}$ 
which is not graded commutative. 
Note that $A_L \cong \Omega^0 \otimes \gl(V)$ as Lie algebras. 

\Proposition A.3. For each $k \in \N$, we obtain cocycles 
$\psi_k \in Z^k(A_L, \Omega^k)$ and 
$\oline\psi_k \in Z^k(A_L, \Omega^{k-1}/d\Omega^{k-2})$ 
satisfying 
$d \circ \oline \psi_k = \psi_k$ by 
$$ \psi_k(f_1 \otimes x_1, \ldots, f_k \otimes x_k) 
:= \sum_{\sigma \in S_k} \Tr(x_{\sigma(1)} \cdots x_{\sigma(k)}) df_1 \wedge \ldots \wedge df_k $$
and 
$$ \oline\psi_k(f_1 \otimes x_1, \ldots, f_k \otimes x_k) 
:= \sum_{\sigma \in S_k} \Tr(x_{\sigma(1)} \cdots x_{\sigma(k)}) [f_1 \cdot 
df_2 \wedge \ldots \wedge df_k] \in \Omega^{k-1} /d\Omega^{k-2}.  $$

\Proof. We first consider the linear map 
$$ \phi = d_B\: A_L \to \Omega^1  \otimes \End(V), \quad 
f \otimes x \mapsto df \otimes x  $$
and note that 
$$(d_{A_L}\phi)(a,b) = a.db - b.da - d[a,b] = [a,db] - [b,da] - [da,b] - [a,db] = 0 $$
implies that $\phi$ is a $1$-cocycle. This implies that its 
$\wedge$-powers $\phi^k \in Z^k(A_L, \Omega^k  \otimes \End(V))$ are 
$k$-cocycles and hence that 
$$ \psi_k = \Tr \circ \phi^k  \in Z^k(A_L, \Omega^k ), $$
because $\Tr \: \Omega^k  \otimes \End(V) \to \Omega^k $ is a morphism 
of the $A_L$-module $\Omega^k  \otimes \End(V)$ onto the trivial module $\Omega^k$. 

That $\oline \psi_k$ is alternating follows from 
$$ [f_1 \wedge df_2 \wedge \ldots \wedge df_k] + [f_2 \wedge df_1 \wedge df_3 \wedge 
\ldots \wedge df_k] 
= [d(f_1f_2)\wedge \ldots \wedge df_k] = [d(f_1f_2 \wedge df_3 \wedge \ldots df_k)] = 0. $$
It remains to show that $\oline \psi_k$ is a $k$-cocycle. To this end, we consider 
$$ \id_A \wedge \phi^{k-1} \in C^{k}(A_L, \Omega^{k-1} \otimes \End(V)) $$
and observe that 
$$ d_{A_L}(\id_A \wedge \phi^{k-1}) = d_{A_L}(\id_A) \wedge \phi^{k-1} = b_A \wedge \phi^{k-1} 
= \id_A \wedge \id_A \wedge \phi^{k-1}. $$
If $q \: \Omega^k  \to \Omega^k /d\Omega^{k-1} , \beta \mapsto [\beta]$ is the quotient map, we have 
$$ \eqalign{ 
&\ \ \ \ \big(q \circ \Tr \circ (\id_A \wedge \phi^{k-1})\big)(f_1 \otimes x_1, \ldots, f_k \otimes x_k) \cr
&= \sum_{\sigma \in S_{k}} \sgn(\sigma) \Tr(x_{\sigma(1)}\cdots x_{\sigma(k)}) 
[f_{\sigma(1)} d f_{\sigma(2)} \wedge \ldots \wedge df_{\sigma(k)}] \cr
&= \sum_{\sigma \in S_{k}} \Tr(x_{\sigma(1)}\cdots x_{\sigma(k)}) 
[f_1 d f_2 \wedge \ldots \wedge df_k] 
= \oline\psi_{k}(f_1 \otimes x_1, \ldots, f_k \otimes x_k). \cr} $$
To see that $\oline \psi_{k}$ is a cocycle, it therefore suffices to show that 
$$ d\oline\psi_k = q \circ \Tr \circ (\id_A \wedge \id_A \wedge \phi^{k-1}) = 0. $$
Explicitly, we have 
$$ \eqalign{ 
&\ \ \ \ \big(q \circ \Tr \circ (\id_A^2 \wedge \phi^{k-1})\big)(f_1 \otimes x_1, \ldots, 
f_{k+1} \otimes x_{k+1}) \cr
&= \sum_{\sigma \in S_{k+1}} \sgn(\sigma) \Tr(x_{\sigma(1)}\cdots x_{\sigma(k+1)}) 
[f_{\sigma(1)} f_{\sigma(2)} \ df_{\sigma(3)}\wedge \ldots \wedge df_{\sigma(k+1)}]. \cr} $$
If $\tau = (1\ 2\ \ldots\ k+1)$ is the cyclic permutation, we obtain 
$$ \eqalign{ 
&\ \ \ \ \sgn(\sigma\tau) \Tr(x_{\sigma\tau(1)}\cdots x_{\sigma\tau(k+1)}) 
[f_{\sigma\tau(1)} f_{\sigma\tau(2)} \ df_{\sigma\tau(3)}\wedge \ldots \wedge df_{\sigma\tau(k+1)}] \cr 
&= (-1)^k\sgn(\sigma)\Tr(x_{\sigma(1)}\cdots x_{\sigma(k+1)}) 
[f_{\sigma(2)} f_{\sigma(3)} \ df_{\sigma(4)}\wedge \ldots \wedge df_{\sigma(k+1)} \wedge df_1] \cr 
&= \sgn(\sigma)\Tr(x_{\sigma(1)}\cdots x_{\sigma(k+1)}) 
[f_{\sigma(2)} f_{\sigma(3)} \ df_{\sigma(1)} 
\wedge df_{\sigma(4)}\wedge \ldots \wedge df_{\sigma(k+1)}].\cr} $$
Now the relation 
$$ [f_1 f_2 df_3 \wedge d(\cdots)] + 
[f_2 f_3 df_1 \wedge d(\cdots)] + 
[f_3 f_1 df_2 \wedge d(\cdots)] + 
= [d(f_1f_2f_3) \wedge d(\cdots)] = 
[d(f_1f_2 f_3 \cdots)] = 0 $$
implies that $3 \oline\psi_k$ is a cocycle, so that $\oline\psi_k$ is a cocycle. 
\qed

\sectionheadline{Appendix B. Cohomology of formal vector fields} 

 In this section, we will review results of Gelfand--Fuks 
on the cohomology of formal vector fields because this is needed to reduce 
Theorem~III.1 to Tsujishita's result. 

Fix $N \in \N$ and write $\F_N = \R[[x_1,\ldots, x_N]]$  for the 
commutative algebra of formal power series in the variables $x_1, \ldots x_N$, endowed 
with the projective limit topology, and let $W_N$ be the Lie algebra of continuous 
derivations of $\F_N$:
$$ W_N = \bigoplus\limits_{j=1}^N \F_N {\d \over \d x_j} .$$
Consider also the subalgebra $L_0$ in $W_N$ consisting of vector fields that vanish
at the origin:
$$ L_0 = \Big\{ \sum_{j = 1}^N f_j (x) {\d \over \d x_j} \quad \: \quad
f_j(0) = 0, j=1, \ldots, N \Big\},  
\quad \hbox{ satisfying } \quad 
W_N \cong L_0 \rtimes \R^N,  $$
where the second factor corresponds to the abelian Lie subalgebra of constant (formal) vector fields. 

For $X \in W_N$, we write 
$$J(X) \in \F_N \otimes \gl_N(\R)$$ 
for the {\it Jacobian of $X$}, defined for 
$X = \sum_i f_i {\d \over \d x_i}$ by $J(X)_{ij} = {\partial f_i \over \partial x_j}$. 
Then 
$$ -J \: W_N \to \F_N \otimes \gl_N(\R)$$ is a crossed homomorphism 
whose kernel is the subalgebra of constant vector fields. 
From the relation 
$$ -J([X,Y])= -X.J(Y) + Y.J(X) + [-J(X),-J(Y)], $$
it follows that 
$$ -J_0 \: L_0 \to \gl_N(\R), \quad X \mapsto - J(X)(0) $$
is a surjective homomorphism of Lie algebras, restricting to an isomorphism on the subalgebra 
of linear vector fields. 

Let $\Omega^p_N$ be the space of formal $p$-forms in $N$ variables, considered as a 
free $\F_N$-module and also as a $W_N$-module, with respect to the action defined by the 
Lie derivative.  

\Definition B.1.  Using the crossed homomorphism $\theta := -J$, 
we get with Lemma~A.1, Proposition~A.3 and Theorem I.7 the following cocycles 
$$ \Phi_k^W \in Z^{2k-1}_c (W_N, \F_N), \quad 
\Phi_k^W (X_1, \ldots, X_{2k-1}) = \sum_{\sigma\in S_{2k-1}} \sgn(\sigma)
\Tr \left( \theta(X_{\sigma(1)}) \ldots \theta(X_{\sigma(2k-1)}) \right) ,$$
$$ \Psi_k^W \in Z^k_c (W_N, \Omega^k_N), \quad 
\Psi_k^W (X_1, \ldots, X_k) = \sum_{\sigma\in S_k} \sgn(\sigma)
\Tr \left( d\theta(X_{\sigma(1)}) \wedge \ldots \wedge d\theta(X_{\sigma(k)}) \right) ,$$ 
and $\Bp_k^W \in Z^k_c (W_N, \Omega^{k-1}_N / d \Omega^{k-2}_N )$, defined by 
$$ \Bp_k^W (X_1, \ldots, X_k) = \sum_{\sigma\in S_k} \sgn(\sigma)
[\Tr \left( \theta(X_{\sigma(1)}) d\theta(X_{\sigma(2)}) \wedge \ldots \wedge d\theta(X_{\sigma(k)}) \right)].
\qeddis 

The following theorem is a re-formulation of the results of Gelfand-Fuks ([GF70a]), 
describing the cohomology of $W_N$ with values in the modules $\Omega_N^p$.
Let $V = \R^N$, considered as the canonical module of $\gl_N(\R)$ and write 
$$\ev_0 \: \Omega_N^p \to \Lambda^p(V'),\quad \omega \mapsto \omega(0) $$
for the evaluation map. We consider $V$ as a module of $L_0$ by pulling back the module 
structure from $\gl_N(\R)$, so that $\ev_0$ is a morphism of $L_0$-modules.  

\Theorem B.2. 
\litem{(a)} For each $p,q \in \N_0$, the map 
$C^p_c(W_N,\Omega^q_N) \to C^p_c(L_0, \Lambda^q(V')), \omega \mapsto  
\ev_0 \circ \omega\res_{L_0}$
induces an isomorphism 
$$ H^p_c(W_N,\Omega^q_N) \to H^p_c(L_0, \Lambda^q(V')). $$
\litem{(b)} $H_c^\bullet (L_0, \Lambda^\bullet (V^\prime)) \cong H_c^\bullet (L_0,\R) 
\otimes H_c^\bullet (L_0, \gl_N(\R), \Lambda^\bullet (V^\prime))$ as bigraded algebras.
\litem{(c)} The inclusions $\gl_N(\R) \into L_0 \into W_N$ induce isomorphisms 
of graded algebras 
$$H_c^\bullet (W_N,\R) \to H_c^\bullet(L_0,\R) \to  H_c^\bullet(\gl_N(\R),\R). $$ 
In particular, $H_c^\bullet (W_N,\R)$ is 
an exterior algebra with generators of degree $2k-1$, $k =1,\ldots, N$. 
\litem{(d)} Let $\Psi_k^L := \ev_0 \circ \Psi_k^W \res_{L_0}
\in Z^k_c (L_0, \gl_N(\R), \Lambda^k (V^\prime))$. 
Then $H_c^\bullet (L_0, \gl_N(\R), \Lambda^\bullet (V^\prime))$ is a quotient 
of the commutative algebra generated by the cohomology classes $[\Psi_k^L]$, 
$k =1,\ldots, N$, by the ideal spanned by the elements of degree exceeding $N$.

\Proof. We only explain how this can be derived from [Fu86]. 

(a) Write $U(\g)$ for the enveloping algebra of a Lie algebra $\g$. 
First we observe that the map 
$$ \eqalign{ \zeta &\: \Omega^p_N \to \Hom_{L_0}(U(W_N), \Lambda^p(V')) \cong \Hom(S(\R^N), \Lambda^p(V')) 
\cong {\cal F}_N \otimes \Lambda^p(V'), \cr 
&\zeta(\alpha)(D) := (D.\alpha)(0) \cr} $$
is an isomorphism of $W_N$-modules. 
Hence $\Omega^p_N$ is coinduced, as a  $W_N$-module, from the 
$L_0$-module $\Lambda^p(V')$. Note that, since $L_0$ is of finite codimensional in $W_N$, 
no problems arise from continuity requirements. 
Therefore (a) follows from a general result on coinduced modules 
([Fu86, Th.~1.5.4]). 

(b) The proof of part (b) is based on the Hochschild--Serre Spectral Sequence 
associated with the
filtration on $C^\bullet (L_0, \Lambda^\bullet (V^\prime))$ relative to the subalgebra $\gl_N(\R) \subset L_0$.
By Theorem 1.5.1(ii) in [Fu86], 
$H^p (L_0, \Lambda^\bullet (V^\prime))$ is the $E_2^{p,0}$-term 
in the spectral sequence, while in the proof of Theorem 2.2.$7^\prime$ in [Fu86], it is shown
that $E_2^{p,0} = E_1^{p,0}$, and the term $E_1^{p,0}$ is calculated explicitly.

(c), (d) follow from (a), combined with [Fu86, Thms.~2.2.7 and 2.2.7']. 

From the formulas in [Fu86] it is not completely obvious that our formula for 
$\Psi_k^L$ 
describes the same cocycle (up to a scalar factor). Fuks describes $\Psi_k^L$ as 
an element of 
$$\Lambda^k(V') \otimes \Lambda^k(S^2(V) \otimes V').$$ 
In our context, 
the cocycle $\Psi_k^L \in Z^k_c(L_0, \gl_N(\R), \Lambda^k (V^\prime))$ is given by the formula 
$$ \Psi_k^L (X_1, \ldots, X_k) = (-1)^k \sum_{\sigma\in S_k} \sgn(\sigma)
\Tr \left( dJ(X_{\sigma(1)})(0) \wedge \ldots \wedge dJ(X_{\sigma(k)})(0)\right) 
\in \Lambda^k(V'). $$
For each $X \in L_0$, the constant term of the 
$1$-form $d J(X) \in \Omega^1_N \otimes \gl_N(\R)$ corresponds to the 
quadratic term in $X$, which can be identified with an element of 
$S^2(V') \otimes V \cong \Sym^2(V) \otimes V$. 

For the basis element $X_{i_1,i_2,\ell} := x_{i_1} x_{i_2} \partial_\ell$, 
we have, in terms of the matrix basis $E_{ij}$ of $\gl_N(\R)$: 
$$ J(X) = x_{i_1} E_{\ell i_2} + x_{i_2} E_{\ell i_1}
\quad \hbox{ and } \quad 
d(J(X)) = dx_{i_1} E_{\ell i_2} + dx_{i_2} E_{\ell i_1} 
= \sum_{\tau \in S_2} dx_{i_{\tau(1)}} E_{\ell i_{{\tau(2)}}}. $$
This leads to 
$$ \eqalign{ 
&\ \ \ \ (-1)^k \Psi_k^L(X_{i_{11},i_{12}, \ell_1}, \ldots, X_{i_{k1}, i_{k2}, \ell_k}) \cr 
&= \sum_{\sigma \in S_k} \sum_{\tau_i \in S_2\atop i=1,\ldots, k} 
 \Tr \left( dx_{i_{\sigma(1),\tau_1(1)}} E_{\ell_{\sigma(1)} i_{\sigma(1)\tau_1(2)}} 
\wedge \cdots \wedge 
dx_{i_{\sigma(k),\tau_k(1)}} E_{\ell_{\sigma(k)} i_{\sigma(k)\tau_k(2)}}\right) \cr
&= \sum_{\sigma \in S_k} \sum_{{\tau_i \in S_2}\atop {i=1,\ldots, k}} 
\delta_{i_{\sigma(1)\tau_1(2)}, \ell_{\sigma(2)}}
\delta_{i_{\sigma(2)\tau_2(2)}, \ell_{\sigma(3)}}\cdots 
\delta_{i_{\sigma(k)\tau_k(2)}, \ell_{\sigma(1)}}\cdot 
dx_{i_{\sigma(1),\tau_1(1)}} \wedge \cdots \wedge dx_{i_{\sigma(k),\tau_k(1)}}. \cr} $$
From this formula it is not hard to verify that our $\Psi_k^L$ are multiples of those in [Fu86]. 
\qed

\sectionheadline{Appendix C. Higher frame bundles and differential forms} 

In this appendix we shall prove our key Theorem III.1. The major part of the 
work consists in explaining why Theorem~III.1 can be derived from 
Tsujishita's [Tsu81]. For that we need the passage from the 
cohomology of formal vector fields to the cohomology of 
${\cal V}_M$ and finally to link it to Theorem~III.1. 
\msk 

Let $M$ be an $N$-dimensional smooth manifold. 
For $k \in \N_0 \cup \{\infty\}$ we write 
$J^k(M)$ for its $k$-frame bundle whose elements are $k$-jets 
$[\alpha]$ of local diffeomorphism $\alpha \: (\R^N,0) \to M$. 
Then $J^0(M) = M$ and $J^1(M)$ is the frame bundle of $M$ whose elements 
are the isomorphisms $\R^N \to T_m(M)$, $m \in M$. 
Evaluating $[\alpha]$ in $0$ leads to a natural map $J^k(M) \to M$ which exhibits 
$J^k(M)$ as a fiber bundle over $M$. 
We write $J^k_m(M)$ for the fiber over $m \in M$. 
For $\ell \leq k$ we have projection maps 
$$ \pi_{\ell,k} \: J^k(M) \to J^\ell(M), \quad [\alpha] \mapsto j_0^\ell(\alpha). $$
In particular, we have 
$\pi_{0,k}([\alpha]) = \alpha(m)$ and 
$\pi_{1,k}([\alpha]) = T_0(\alpha)$. 

\msk
Let $G^k$ denote the group of $k$-jets of local diffeomorphisms of $\R^N$ fixing $0$. 
As a set, $G^k$ is the set of all polynomial maps 
$\phi \: \R^N \to \R^N$ of degree $\leq k$ without constant term for which 
$d\phi(0)$ is invertible. The multiplication is given by 
$$ \phi_1\phi_2 := j^k_0(\phi_1 \circ \phi_2). $$
For $\ell \leq k$ we have natural surjective homomorphisms 
$$ q_{\ell,k} \: G^k \to G^\ell, \quad \alpha \mapsto j^\ell_0(\alpha), $$
cutting off all terms of order $> \ell$. We put 
$$ G_k := \ker(q_{k,\infty}) \trile G^\infty $$
and note that 
$$ G^\infty \cong G_1 \rtimes \GL_N(\R) \leqno(C.1) $$
in a natural way. 

\msk 
For $k < \infty$, the group $G^k$ is a finite-dimensional Lie group, 
$G^1 \cong \GL_N(\R)$, and 
$$ G^k \cong G^k_1 \rtimes \GL_N(\R), $$
where $G^k_1 := \ker(q_{1,k})$ is a simply connected nilpotent Lie group. 
Identifying $G^\infty$ with the projective limit of all groups 
$G^k \cong G^\infty/G_k$, we obtain a 
topology for which it actually is a Lie group with Lie algebra 
$\Lie(G^\infty) \cong L_0$ and 
$L_k := \Lie(G_k)$ is a finite-codimensional ideal of $L_0$ (cf.\ Appendix B). 
The normal subgroup $G_1$ is pro-nilpotent, its exponential 
function $\exp_{G_1} \: L_1 \to G_1$ is a diffeomorphism, and (C.1) 
is a semidirect product of Lie groups. 

The group  $G^k$ acts on $J^k(M)$ from the right by 
$[\alpha].\phi := [\alpha \circ \phi]$. This action is transitive on the fibers and 
defines on $J^k(M)$ the structure of a smooth $G^k$-principal bundle. 

The group $\Diff(M)$ acts on each frame bundle $J^k(M)$ by 
$\phi.[\alpha] := [\phi  \circ \alpha]$, and the corresponding 
homomorphism of Lie algebras is given by 
$$ \gamma^k \: {\cal V}_M \to {\cal V}_{J^k(M)}, \quad 
\gamma^k(X)_{[\alpha]} = \derat0 [\Fl^X_t \circ \alpha]. \leqno(C.2) $$
Since the action of $\Diff(M)$ on $J^k(M)$ 
commutes with the action of the structure 
group $G^k$, 
$$ \gamma^k({\cal V}_M) \subeq {\cal V}_{J^k(M)}^{G^k}. \leqno(C.3)$$

\Lemma C.1. {\rm([Tsu81, Lemma~4.2.2])} Let $k \in \N \cup \{\infty\}$, 
$\nabla$ be an affine connection on $M$, and 
$D \subeq TM$ the open domain of the 
corresponding exponential function 
$\Exp^\nabla \: D \to M$, which is an open neighborhood of the zero 
section. Then 
$$ s \: J^1(M) \to J^k(M), \quad s(\alpha) := j_0^k(\Exp^\nabla_m \circ \alpha) $$
is a smooth section. 
\qed

Note that for each $m \in M$ the intersection $D_m := D \cap T_m(M)$ is an open 
zero neighborhood and that $\Exp^\nabla_m \: D_m \to M$ is a smooth map 
with $T_0(\Exp^\nabla_m) = \id_{T_m(M)}$, hence a local diffeomorphism. 
Now the map $J^1(M) \times G_1^k \to J^k(M), (\alpha, g) \mapsto s(\alpha)g$ 
is a diffeomorphism and the map 
$$ F \: J^k(M) \to G_1^k, \quad s(\alpha)g \mapsto g, \quad \alpha \in J^1(M), $$
is smooth and $G_1^k$-equivariant. 
For $g \in \GL_N(\R)$ we have 
$s(\alpha \circ g) = s(\alpha).g$, which implies that 
$$ F(\alpha.g) = g^{-1} F(\alpha) g \quad \hbox{ for all } \quad \alpha \in J^k(M). \leqno(C.4)$$
Hence $F$ is equivariant with respect to the action of 
$G^k\cong G_1^k \rtimes \GL_N(\R)$ on $G_1^k$ from the right by 
$g.(g_1, g_0) = g_0^{-1}gg_1 g_0$.

Since Tsujishita uses a realization of Lie algebra cohomology in terms of 
right invariant, resp., equivariant differential forms on the corresponding 
group, we briefly discuss the relevant identifications in the following remark. 

\Remark C.2. (a) Let $G$ be a Lie group. A {\it smooth $G$-module} is a topological 
vector space $V$ on which $G$ acts smoothly by linear maps. We write 
$$ \rho_V \: G \times V \to V, \quad (g,v) \mapsto \rho_V(g)(v) =: g.v$$ 
for the action map. 

Further, let $M$ be a smooth  manifold on which $G$ acts from the right 
by $M \times G \to M, (m,g) \mapsto m.g =: \rho_g^M(m)$. 
We call a $p$-form $\alpha \in \Omega^p(M,V)$ {\it equivariant} if we have 
for each $g \in G$ the relation 
$$ (\rho_g^M)^*\alpha = \rho_V(g)^{-1} \circ \alpha. $$
We write $\Omega^p(M,V)^G$ for the space of $G$-equivariant $p$-forms on $M$. 
This is the space of $G$-fixed elements with 
respect to the action of $G$ on $\Omega^p(M,V)$, 
given by $g.\alpha := \rho_V(g)\circ (\rho_g^M)^*\alpha$. 

(b) For the right action of $G$ on $M = G$ by left multiplication 
$x.g := g^{-1}x$ we obtain the space 
$\Omega^p_\ell(G,V)$ of left equivariant forms, i.e., forms satisfying 
$$ \lambda_g^*\omega = \rho_V(g) \circ \omega, \quad g \in G. $$
In [ChE48] it is shown that the evaluation map 
$$ \ev_\1 \: (\Omega^\bullet_\ell(G,V),d) \to (C^\bullet_c(\g,V),d_\g), \quad 
\omega \mapsto \omega_\1 $$
yields an isomorphism of cochain complexes 
(cf.\ [Ne04] for the unproblematic 
extension to infinite-dimensional Lie groups). 

(c) There is also a realization of the complex $(C^\bullet_c(\g,V),d_\g)$  
by right equivariant differential forms on $G$: If $\eta_G \: G \to G, g \mapsto g^{-1}$ is the 
inversion map, then for each left equivariant $p$-form $\alpha \in \Omega^p_\ell(G,V)$ 
the $p$-form $\tilde\alpha := \eta_G^*\alpha$ is {\it right equivariant}, 
i.e., 
$$ \rho_g^*\tilde\alpha = \rho_V(g)^{-1} \circ \tilde\alpha \quad \hbox{ for each } \quad 
g \in G. $$
We thus obtain an isomorphism of cochain complexes 
$$\eta_G^* \: (\Omega^\bullet_\ell(G,V),d) 
\to (\Omega^\bullet_r(G,V),d). $$
Since $T_\1(\eta_G) = - \id_G$, we also obtain an isomorphism 
$$ \tilde\ev_\1 \: (\Omega^\bullet_r(G,V),d) \to (C^\bullet_c(\g,V),d_\g),
 \quad \omega \mapsto (-\id_\g)^*\omega_\1. $$

For each $g \in G$ and $\omega \in \Omega^p_r(G,V)$ the form 
$\lambda_g^*\omega$ is also right equivariant and satisfies 
$$ (\lambda_{g^{-1}}^*\omega)_\1(x_1,\ldots, x_p) 
= \omega_{g^{-1}}(g^{-1}.x_1,\ldots, g^{-1}.x_p) 
= \rho_V(g).\omega(\Ad(g)^{-1}.x_1,\ldots, \Ad(g)^{-1}.x_p), $$
showing that $\tilde\ev_\1$ intertwines the action of $G$ by left translation 
on $\Omega^p_r(G,V)$ with the natural action of $G$ on the cochain space $C^p_c(\g,V)$. 

(d) There is an alternative realization of the complex $(C^\bullet(\g,V), d_\g)$ 
as the space of left, resp., right invariant $V$-valued differential forms on $G$. 

First we write the Chevalley--Eilenberg differential as 
$$ d_\g = d_\g^0 + \rho_V \wedge, $$
where $d_\g^0$ is the differential corresponding to the trivial module structure on 
$V$ and $\rho_V \: \g \to \gl(V)$ is the homomorphism of Lie algebras, defining the 
$\g$-module structure on $V$. 

If $\kappa_G \in \Omega^1(G,\g)$ denotes the left Maurer--Cartan form of $G$, 
we obtain $\kappa_V := \rho_V \circ \kappa_G \in \Omega^1(G,\gl(V))$ 
and a corresponding covariant differential on $\Omega^\bullet(M,V)$: 
$$ d_\kappa \omega = d\omega + \kappa_V \wedge \omega. $$
Since $\kappa$ is left invariant, this differential preserves the subspace of 
left invariant forms, and it is easy to see that evaluation in $\1$ 
intertwines it with $d_\g$ on $C^\bullet_c(\g,V)$. 

In a similar fashion, one obtains a realization by right invariant differential forms 
with the appropriate differential. 
\qed

Let $A$ be a finite-dimensional smooth $G^\infty$-module, where we assume that 
the action of the non-connected Lie group $\GL_N(\R)$ on $A$ is the one 
obtained by restricting the action of $\GL_N(\C)$ on the complexification $A_\C$, 
hence completely determined by the action of $\gl_N(\R)$ on $A$. 
We then have natural identifications  
$$ C^\bullet_c(L_0, \gl_N(\R), A) \cong 
C^\bullet_c(L_1, A)^{\gl_N(\R)} \cong 
\Omega_r^\bullet(G_1, A)^{\GL_N(\R)} \leqno(C.5) $$
([Tsu81, Lemma~3.3.4], Remark~C1b(c)). 
Although $\GL_N(\R)$ is not connected, this follows 
from a complexification argument since $\GL_N(\C)$ is connected. 

Further, let $\omega \in 
C^p_c(L_1, A)^{\GL_N(\R)} \cong C^p_c(L_1, A)^{\gl_N(\R)}$ and 
$\tilde\omega \in \Omega_r^p(G_1, A)$ the corresponding 
right equivariant $A$-valued $p$-form with 
$\tilde\omega_\1 = (-1)^p \omega$ (cf.\ Remark~C.2(c)). 
For $g_1 \in G_1$ we then have 
$$ \rho_{g_1}^*F^*\tilde\omega 
= (F \circ \rho_{g_1})^*\tilde\omega 
= (\rho_{g_1} \circ F)^*\tilde\omega 
= \rho_A({g_1}) \circ F^*\tilde\omega $$
and for $g_0 \in \GL_N(\R)$ we obtain with (C.4) 
$$ \rho_{g_0}^*F^*\tilde\omega 
= (F \circ \rho_{g_0})^*\tilde\omega 
= (c_{g_0}^{-1} \circ F)^*\tilde\omega 
= F^*\big( (c_{g_0}^{-1})^*\tilde\omega\big)
= \rho_A(g_0) \circ F^*\tilde\omega, $$
hence the $\GL_N(\R)$-equivariance of 
$\tilde\omega$. We thus obtain for $k = \infty$ a morphism of chain complexes 
$$ \tilde F^* \: C^\bullet_c(L_1,A)^{\GL_N(\R)} \to 
\Omega^\bullet(J^\infty(M), A)^{G^\infty}, \quad  
\omega \mapsto F^*\tilde\omega. \leqno(C.6) $$

The map $\gamma_A$ defined below is the one used by Tsujishita in [Tsu81, p.~62]. 

\Proposition C.3. Let ${\cal A} := J^\infty(M) \times_{G^\infty} A \to M$ 
denote the bundle associated to $J^\infty(M)$ by the representation 
of $G^\infty$ on $A$ and $\Gamma({\cal A})$ be its space of smooth sections on which 
${\cal V}_M$ acts via $\gamma^\infty$. We identify $\Gamma({\cal A})$ 
with $C^\infty(J^\infty(M), A)^{G^\infty}$. Then the map 
$$ \gamma_A \:= (\gamma^\infty)^* \circ \tilde F^* \: 
C^\bullet_c(L_1,A)^{\GL_N(\R)} 
\cong C^\bullet_c(L_0,\gl_N(\R),A) 
\to C^\bullet_c({\cal V}_M, \Gamma({\cal A})) 
\leqno(C.7) $$
is a morphism of chain complexes. Different choices of affine connections on 
$M$ yield the same maps 
$$ H^\bullet_c(L_1,A)^{\GL_N(\R)} \cong H^\bullet_c(L_0,\gl_N(\R),A)\to 
H^\bullet_c({\cal V}_M, \Gamma({\cal A}))$$
in cohomology. 

\Proof. For the continuity of the so obtained cochains $\gamma_A(\omega)$ of ${\cal V}_M$, 
we note that the continuous representation of $G^\infty$ on the finite-dimensional space $A$ 
factors through some finite-dimensional 
quotient group $G^k$. Hence the bundle ${\cal A}$ is also 
associated to $J^k(M)$, $\Gamma({\cal A})$ can be realized as $A$-valued 
$G^k$-equivariant functions on $J^k(M)$, and the action of ${\cal V}_M$ on 
this space via $\gamma^k \: {\cal V}_M \to {\cal V}_{J^k(M)}$ is continuous 
because the finite dimensionality of $J^k(M)$ implies that 
$\gamma^k$ is a continuous morphism of topological Lie algebras.

To see that the choice of connection has no effect on the corresponding map in 
cohomology, let $\nabla'$ be another affine connection on $M$ and observe that 
$\nabla_t := t \nabla' + (1-t)\nabla$ defines a smooth family of affine connections 
on $M$ with $\nabla_0 = \nabla$ and $\nabla_1 = \nabla'$. 
It is easy to see that the corresponding functions $s_1 \: J^1(M) \to J^\infty(M)$ 
depend smoothly on $t$, and so do the corresponding $G_1$-equivariant 
functions $F_t \: J^\infty(M) \to  G_1$. Hence, for 
each $\omega \in C^p_c(L_1,A)^{\GL_N(\R)}$ the differential forms 
$F_1^*\tilde\omega - F_0^*\tilde\omega$
are equivariantly exact, i.e., the differential of a $G^\infty$-equivariant 
$A$-valued $(p-1)$-form (cf.\ [Ko74, p.143]). 
This implies the assertion. 
\qed

\Remark C.4. We also note that the $1$-form 
$\delta(F) \in \Omega^1(J^\infty(M),L_1)$
satisfies the Maurer--Cartan equation, hence defines a crossed homomorphism 
$$\delta(F) \: {\cal V}(J^\infty(M)) \to C^\infty(J^\infty(M), L_1)$$ 
(Proposition~I.8). Composing with the homomorphism 
$\gamma^\infty \: {\cal V}_M \to {\cal V}_{J^\infty(M)}$, we 
thus obtain a crossed homomorphism 
$$\delta(F) \circ \gamma^\infty \: {\cal V}_M\to C^\infty(J^\infty(M),L_1)^{G^\infty}.
$$ 
That $\delta(F) \circ \gamma^\infty$ 
maps into $G^\infty$-equivariant functions is due to (C.3) and the 
$G^\infty$-equivariance of $F$, resulting from (C.4). 
We further note that any element $\omega \in C^p_c(L_1,A)^{\GL_N(\R)}$ 
defines by composition a ${\cal V}_M$-equivariant $p$-cochain 
of the Lie algebra $C^\infty(J^\infty(M),L_1)^{G^\infty}$ with values 
in $C^\infty(J^\infty(M), A)^{G^\infty} \cong \Gamma({\cal A})$. 
Now Theorem~I.7 shows that pulling back with the crossed homomorphism 
$\delta(F) \circ \gamma^\infty$ yields $\Gamma({\cal A})$-valued Lie algebra 
cochains and that this is compatible with the mutual differentials. 
\qed

\Example C.5. For $A = \Lambda^k(V')$, $V = \R^N$, and the canonical representation 
of $G^1 \cong \GL_N(\R) \cong G^\infty/G_1$ on this space, we obtain for the space of smooth 
sections $\Gamma({\cal A}) \cong \Omega^k_M,$
so that we get a morphism of cochain complexes 
$$ \gamma_A  \: C^\bullet_c(L_1,\Lambda^k(V'))^{\GL_N(\R)} \to C^\bullet_c({\cal V}_M, \Omega^k_M). $$

In this case the bundle ${\cal A} = \Lambda^k(T^*(M))$ is associated to the 
frame bundle $J^1(M)$ which creates a simpler picture than working with 
the infinite-dimensional manifold $J^\infty(M)$. 
\qed

\Example C.6. If $A = \R$ is the trivial module, we obtain in particular 
${\cal F}_M$-valued cocycles on ${\cal V}_M$ from any map 
$$ (\gamma^k)^* \: (\Omega^p_{J^k(M)})^{G^k} \to C^p_c({\cal V}(M),{\cal F}_M). $$

Here is a concrete example: For $k = 1$ we consider the 
$1$-form $\omega \in \Omega^1(J^1(M),\R)^{\GL_N(\R)}$ defined as follows. 
From the homomorphism 
$$ \chi \: \GL_N(\R) \to \R^\times_+, \quad g \mapsto |\det(g)| $$ 
we obtain an associated bundle 
$J^1(M) \times_\chi \R^\times_+$, and since $\R^\times_+$ is contractible, 
this bundle has a global section (we could also take $\log \circ \chi$ and obtain 
an affine bundle), which means that there is a smooth function 
$$ F \: J^1(M) \to \R^\times_+ \quad \hbox{ with } \quad 
F(\alpha.g) = F(\alpha) |\det(g)| \quad \hbox{ for } 
\quad g \in \GL_N(\R), \alpha \in J^1(M). $$
Then 
$$\delta(F) \in \Omega^1(J^1(M), \R)^{\GL_N(\R)} $$
is a $\GL_N(\R)$-invariant $1$-form. 

If $\mu$ is a volume form on $M$, then we can construct $F$ directly from 
$\mu$ by 
$$ F(v_1,\ldots, v_N) := |\mu_m(v_1,\ldots, v_N)| $$
for any basis $(v_1,\ldots, v_N)$ on $T_m(M)$. 

For a diffeomorphism $\phi \in \Diff(M)$, we then have 
$$ \eqalign{ (\phi.F)(v_1,\ldots, v_N) 
&= F(T(\phi^{-1}).(v_1,\ldots, v_N))
= |\mu(T(\phi^{-1}).(v_1,\ldots, v_N))| \cr
&= |(\phi^{-1})^*\mu(v_1,\ldots, v_N)|. \cr} $$
Dividing by $F(v_1,\ldots, v_N)$, we obtain a smooth function 
$$ (\phi.F)F^{-1} \in C^\infty(M,\R^\times_+), $$
and passing to the Lie derivative, we obtain for each vector field $X$ on $M$ 
a smooth function 
$$ ({\cal L}_X.F)F^{-1} \in C^\infty(M,\R). 
\qeddis

Now we turn to Tsujishita's construction of the homomorphism 
$$ H^\bullet_c(L^1,\Lambda^k(V'))^{\GL(V)} \to H^\bullet_c({\cal V}_M, \Omega^k_M) $$
in term of the cocycles ${\cal L}_X\nabla$ associated to affine 
connections (Lemma~II.1). 

Let $\nabla$ be an affine connection on $M$. As in Lemma~C.1, we obtain from $\nabla$ 
a smooth section 
$$ s \: J^1(M) \to J^2(M), \quad s(\alpha) 
:= j_0^2(\Exp^\nabla_m \circ \alpha) = j_0^2(\Exp^\nabla_m) \circ \alpha. $$
From $s(\alpha g) = s(\alpha)g$ for $\beta \in \GL(V)$ it follows that 
$s$ is $\GL(V)$-equivariant. 

Let $F \: J^2(M) \to G_1/G_2 \cong \Sym^2(V,V)$ denote the unique smooth $G^1$-equivariant smooth 
function vanishing on $s(J^1(M))$. Identifying $\Sym^2(V,V) \cong L_1/L_2$ with the corresponding subspace of $V' \otimes \gl(V) \cong (V' \otimes V') \otimes V$ corresponds to composition 
with the map $\ev_0 dJ \: L_1 \to V' \otimes \gl(V)$.  

The map 
$$\delta(F) \circ \gamma^2 \: {\cal V}_M \to C^\infty(J^2(M), L_1/L_2)^{G^2} $$ 
is a crossed homomorphism (Proposition~I.8), and since $L_1/L_2$ is abelian, it is 
a $1$-cocycle. 
Composing with $\ev_0 dJ \: L_1/L_2 \to V' \otimes\gl(V)$, we thus get a $1$-cocycle 
$$\eqalign{ \ev_0 dJ\delta(F) \circ \gamma^2 
&\: {\cal V}_M \to C^\infty(J^2(M), V' \otimes \gl(V))^{G^2} \cr
&\qquad\qquad \cong C^\infty(J^1(M), V' \otimes \gl(V))^{\GL(V)} \cong \Omega^1(M,\End(TM)). \cr}$$

The following theorem is the link between the approaches in [Tsu81] and [Ko74]. 

\Theorem C.7. If $\nabla$ is torsion free, then 
$$ \ev_0 dJ\delta(F)(\gamma^2(X)) = {\cal L}_X \nabla \in \Omega^1(M,\End(TM)). $$

\Proof. Since this is a local assertion, it suffices to assume that 
$M$ is an open subset of $V \cong \R^N$. Then we write 
$$ \nabla_X Y = dY(X) + \Gamma(X,Y) $$
for a symmetric $(1,2)$-tensor $\Gamma$ (recall that $\Gamma$ is symmetric if and only if 
$\nabla$ is torsion free). 
We then identify  
$$ J^1(M) = M \times \GL(V) \quad \hbox{ and } \quad 
J^2(M) = M \times G^2 \cong M \times (\Sym^2(V,V) \rtimes \GL(V)). $$
Here we write elements of $G^2$  as 
$\alpha = \alpha_1 + \alpha_2$ with $\alpha_1 \in \GL(V)$ linear and 
$\alpha_2$ quadratic. Then the group structure of $G^2$ is given by 
$$ \alpha \bullet \beta 
= \alpha_1 \beta_1 + (\alpha_1 \beta_2 + \alpha_2 \beta_1) 
= \alpha_1 \beta + \alpha_2 \beta_1 
= \alpha \beta_1 + \alpha_1 \beta_2. $$
From this we see that 
$$ \alpha^{-1} = \alpha_1^{-1} - \alpha_1^{-1} \alpha_2 \alpha_1^{-1}. \leqno(C.10)$$

To describes the section 
$s \: J^1(M) \to J^2(M)$
explicitly, let $\gamma_v(t) = \Exp^\nabla(tv)$ be the geodesic 
with 
$\gamma_v(0) = m$ and 
$\gamma_v'(0) = v$ (which is defined for $t$ sufficiently close to $0$). 
Then the relation 
$$ \gamma_v''(t) = - \Gamma(\gamma_v'(t),\gamma_v'(t)) $$ 
leads to $\gamma_v''(0) = - \Gamma(v,v)$, so that 
$$ j_0^2(\Exp^\nabla_m)(v) = v + \half\gamma''(0) = v - \half\Gamma(v,v). $$
We thus obtain 
$$ s(m,\alpha_1) = (m, (\1 - \half\Gamma)\alpha_1). $$
Accordingly, we find 
$$ F(m,\alpha) 
= \alpha_1^{-1}(\alpha +\half \Gamma \circ \alpha_1) 
= \1 + \alpha_1^{-1}\alpha_2 + \half \alpha_1^{-1}\Gamma\circ \alpha_1 $$
because 
$$ F(m,\alpha(\1 + \beta_2)) 
=  F(m,\alpha + \alpha_1 \beta_2) 
=  F(m,\alpha) + \beta_2 = F(m,\alpha).\beta_2 $$
and 
$$ F(m,\alpha\beta_1) = \beta_1^{-1} F(m,\alpha)\beta_1. $$
From the relation 
$$ (\beta \bullet \alpha)_1^{-1} (\beta \bullet \alpha)_2 
= \alpha_1^{-1} \beta_1^{-1} (\beta_1 \alpha_2 + \beta_2 \alpha_1) 
= \alpha_1^{-1} \alpha_2 + \alpha_1^{-1} \beta_1^{-1}\beta_2 \alpha_1, $$
we obtain the following formula for the differential of $F$ in $(v,\beta \bullet \alpha) \in 
T_{(m,\alpha)}(M \times G^2)$ with $\beta \in T_\1(G^2) \cong L_0/L_2$: 
$$ \eqalign{ dF(m,\alpha)(v,\beta\bullet \alpha) 
&= \half \alpha_1^{-1} d\Gamma(v) \circ \alpha_1 
+ \alpha_1^{-1} \beta_2 \alpha_1 - \half \alpha_1^{-1} (\beta_1.\Gamma) \circ \alpha_1  \cr
&= \alpha_1^{-1} \Big(\half d\Gamma(v) + \beta_2  - \half (\beta_1.\Gamma)\Big) \circ \alpha_1. \cr}$$  

The natural lift of $X \in {\cal V}_M$ to $J^2(M)$ is given by 
$$ \eqalign{ \gamma^2(X)(m,\alpha) 
&= j_0^2(X \circ \alpha) 
= (X(m), J(X)_m \circ \alpha + \half (d^2X)_m \circ \alpha_1)\cr
&= \Big(X(m), \big(J(X)_m + \half (d^2X)_m\big) \bullet \alpha\Big) \cr} $$
($J(X)_m v = (dX)_m(v)$ denoting the Jacobian of $X$), 
so that 
$$\gamma^2(X)(s(m,\alpha_1))
= \Big(X(m), J(X)_m \circ \alpha_1 - \half J(X)_m \Gamma\circ \alpha_1 
+ \half (d^2X)_m \circ \alpha_1\Big). $$
This leads to 
$$ \eqalign{ 
dF(\gamma^2(X))(s(m,\alpha_1))
&= \delta(F)(\gamma^2(X))(s(m,\alpha_1)) \cr
&= \half \alpha_1^{-1} \Big(d\Gamma(X(m)) + (d^2X)_m - J(X)_m.\Gamma\Big)\circ \alpha_1.\cr} $$
This is a $\GL(V)$-equivariant function $J^1(M) \to \Sym^2(V,V)$, so that the corresponding 
$(1,2)$-tensor field is given in the canonical coordinates by the smooth function 
$M \to \Sym^2(V,V)$ given by 
$$ m \mapsto dF(\gamma^2(X))(s(m,\1))= \half\big(d\Gamma(X(m)) + (d^2X)_m - J(X)_m\Gamma\big). $$

Suppose that $\beta \: V \to V$ is a quadratic map and 
$\tilde\beta$ 
the corresponding symmetric bilinear map determined by 
$\beta(v) = \half\tilde\beta(v,v)$. Considering $\beta$ as an element of $L_0/L_2$, 
we have 
$$ J(\beta)(v) = (d\beta)_v = \tilde\beta(v,\cdot) \in \gl(V) $$ 
and thus 
$$ (dJ(\beta))_0(v)v = \ev_0(dJ(\beta))(v)(v) = \tilde\beta(v,v) = 2 \beta(v). $$
Applying $\ev_0 \circ dJ$ to the smooth function above thus leads to 
$$ \ev_0dJ \delta(F)(\gamma^2(X))(s(m,\1))= d\Gamma(X(m)) + (d^2X)_m - J(X)_m\Gamma \in V' \otimes \gl(V). $$

In view of Remark~II.2, it remains to show that 
$$ {\cal L}_X\Gamma = d\Gamma(X)  - J(X).\Gamma. $$
It suffices to verify this relation by evaluating on 
constant vector fields $Y$ and $Z \in V$: 
$$ \eqalign{ ({\cal L}_X\Gamma)(Y,Z) 
&= [X,\Gamma(Y,Z)] - \Gamma([X,Y],Z) - \Gamma(Y,[X,Z]) \cr
&= (d\Gamma(X))(Y,Z) - dX(\Gamma(Y,Z)) +\Gamma(dX(Y), Z)+ \Gamma(Y,dX(Z))\cr
&= (d\Gamma(X))(Y,Z) - (J(X).\Gamma)(Y,Z). \cr} $$
\qed

Now we can identify the image of the $\Lambda^k(V')$-valued cocycle 
$\Psi_k^L$ (Theorem~B.2) 
under Tsujishita's map $\gamma_{\Lambda^k(V')}$, described in Proposition~C.3,  
which maps it into an element of $Z^k_c({\cal V}_M, \Omega^k_M)$. 

\Theorem C.8. For each torsion free connection $\nabla$ on the manifold $M$, 
we have 
$$ \gamma_{\Lambda^k(V')}(\Psi_k^L) = \Psi_k. $$

\Proof. We first observe that 
$$\Psi_k^L= (\ev_0dJ)^*f_k, \quad \hbox{ where } \quad 
f_k \: (V' \otimes \End(V))^k \to \Lambda^k(V') $$
is given by 
$$ \eqalign{ ((\alpha_1 \otimes A_1), \ldots, (\alpha_k \otimes A_k)) 
&\mapsto 
\sum_{\sigma\in S_k} \sgn(\sigma)
\Tr(A_{\sigma(1)} \cdots A_{\sigma(k)}) \alpha_{\sigma(1)} \wedge \cdots \wedge 
\alpha_{\sigma(k)}  \cr
&= \Big(\sum_{\sigma\in S_k} \Tr (A_{\sigma(1)} \cdots A_{\sigma(k)})\Big) 
\cdot \alpha_1 \wedge \cdots \wedge \alpha_k. \cr} $$
In this sense, we have $f_k =m_k \otimes \beta_k$, where 
$m_k \: (V')^{\otimes k} \to \Lambda^k(V')$ is the alternating multiplication map and 
$\beta_k \: \End(V)^k \to \R$ is the symmetric $k$-linear map defined by 
$$ \beta_k(A_1, \ldots, A_k) = \sum_{\sigma\in S_k} \Tr(A_{\sigma(1)} 
\cdots A_{\sigma(k)}). $$

According to Proposition~C.3, we have for $X_1, \ldots, X_k \in {\cal V}_M$: 
$$ \gamma_{\Lambda^k(V')}(\Psi_k^L)(X_1, \ldots, X_k) 
= (F^*\tilde\Psi_k^L)(\gamma^2(X_1), \ldots, \gamma^2(X_k)), $$
where we identify 
$$\Omega^k_M 
\cong C^\infty(J^\infty(M), \Lambda^k(V'))^{G^\infty} 
\cong C^\infty(J^2(M), \Lambda^k(V'))^{G^2} 
\cong C^\infty(J^1(M), \Lambda^k(V'))^{\GL(V)}. $$
This identification is based on the canonical map 
$$ J^\infty(M) \times \Lambda^k(V') \to \Lambda^k(T^*(M)), \quad 
([\alpha],\omega) \mapsto \omega \circ \alpha_1^{-1} $$
which is constant on the diagonal $G^\infty$-action and factors through 
the corresponding map on $J^1(M)$, resp., $J^2(M)$. 

Since the cocycles $\Psi_k^L$ factor through the quotient algebra $L_0/L_2$, 
the map $\gamma_{\Lambda^k(V')}$ can be constructed with $J^2(M)$ instead of 
$J^\infty(M)$ and to pass from $G^2$-equivariant smooth functions on $J^2(M)$ 
to $\GL(V)$-equivariant smooth functions on $J^1(M)$, we may simply 
restrict to $s(J^1(M))$, the $\1$-level set of the function $F$. We thus obtain with 
Theorem~C.7 for any torsion free connection $\nabla$: 
$$ \eqalign{ 
\gamma_{\Lambda^k(V')}(\Psi_k^L)(X_1, \ldots, X_k) 
&= (F^*\tilde\Psi_k^L)(\gamma^2(X_1), \ldots, \gamma^2(X_k)) \cr
&= (-1)^k \Psi_k^L(\delta(F)(\gamma^2(X_1)), \ldots, 
\delta(F)(\gamma^2(X_k)))\cr
&= (-1)^k f_k\Big(\ev_0dJ\delta(F)(\gamma^2(X_1)), \ldots, 
\ev_0dJ\delta(F)(\gamma^2(X_k))\Big) \cr
&= (-1)^k f_k({\cal L}_{X_1}\nabla, \ldots, {\cal L}_{X_k}\nabla) \cr
&= \Psi_k(X_1,\ldots, X_k), \cr} $$
where we identify ${\cal L}_X\nabla$ with an element of $\Omega^1(M,\End(TM)) 
\cong C^\infty(J^1(M), V' \otimes \End(V))^{\GL(V)}$ and use 
$\Omega^k_M \cong C^\infty(J^1(M), \Lambda^k(V'))^{\GL(V)}$. 
\qed

\nin{\bf Proof of Theorem III.1:} We choose an affine torsion free connection $\nabla$. 
From the formulation of Tsujishita's Theorem in [Tsu81], Thm.~5.1.6, 
we get (with $V = \R^N$) an isomorphism of bigraded algebras: 
$$H^\bullet_c({\cal V}_M, \Omega^\bullet_M) \cong 
H^\bullet_c ({\cal V}_M, \F_M) \otimes H^\bullet_c (L_0, \gl_N(\R), 
\Lambda^\bullet (V^\prime)).$$
We further know from Theorem~B.2, that $H^\bullet_c (L_0, \gl_N(\R), 
\Lambda^\bullet (V^\prime))$ is generated by the classes of the cocycles 
$\Psi_k^L$. Hence it suffices to observe with Theorem~C.8 that Tsujishita's map 
$$ \gamma_{\Lambda^k(V')} \: C^\bullet_c(L_1,\Lambda^k(V'))^{\GL_N(\R)} 
\cong C^\bullet_c(L_0,\gl_N(\R),\Lambda^k(V')) 
\to C^\bullet_c({\cal V}_M, \Omega^k_M) $$ 
(Proposition~C.3) maps $\Psi_k^L$ to $\Psi_k$. 
\qed

\sectionheadline{Appendix D. Cohomology of vector fields with values in the trivial module} 

Here we shall review the results of Haefliger [Hae76]. 

We begin by recalling the definition of the Weil algebra $\WW$.
Let $E(u_1,\ldots, u_N)$ be the exterior algebra in the generators $u_1, \ldots, u_N$, and let 
$\R[c_1, \ldots, c_N]$, be the polynomial algebra in the generators $c_1,\ldots, c_N$. Introduce 
gradings on these algebras by assigning degrees to the generators:
$$ \deg(u_k) = 2k-1 \quad \hbox{ and }  \quad \deg(c_k) = 2k. $$
Consider the quotient ${\widehat \R} [c_1, \ldots, c_N]$ of the polynomial algebra $\R[c_1, \ldots, c_N]$ 
by the ideal spanned by the elements of degrees exceeding $2N$.
The {\it Weil algebra} is the differential graded algebra
$$ \WW = E(u_1,\ldots, u_N) \otimes {\widehat \R} [c_1, \ldots, c_N]$$
with the differential defined by
$$ d(u_k) = c_k \quad \hbox{ and } \quad d(c_k) = 0.$$
 Let $H^\bullet(\WW)$ be the cohomology of this 
differential graded algebra.
 
 Gelfand--Fuks used the Weil algebra to describe the cohomology of the Lie algebra $W_N$ 
of formal vector fields in dimension~$N$:
 
\Theorem D.1. {\rm([GF70b]) $H^\bullet_c(W_N,\R) \cong H^\bullet(\WW)$} 
 
 The explicit description of the cohomology of the Weil algebra was given by Vey:
 
\Theorem D.2. {\rm([Go74])} $H^\bullet (\WW)$ is spanned by $1$ and the cocycles 
$$ (u_{i_1} \wedge \ldots \wedge u_{i_r}) \otimes (c_{j_1} \ldots c_{j_s}) \leqno(D.1) $$
satisfying the following conditions: $1 \leq i_1 < \ldots < i_r \leq N$;
$1 \leq j_1 \leq \ldots \leq j_s \leq N$; $i_1 \leq j_1$; $r > 0$; $j_1 + \ldots + j_s \leq N$;
$i_1 + j_1 + \ldots + j_s > N$.
\qed

Combining the preceding two theorems, we derive: 
\Corollary D.3. 
\litem{(a)} $H^q (W_N,\R) = 0$ for $1 \leq q \leq 2N$ and for $q \geq (N+1)^2$.
\litem{(b)} $\dim H^{2N+1} (W_N,\R) = p(N+1) -1$, where $p$ is the partition function.
\litem{(c)} $\dim H^{2N+2} (W_N,\R) = 0$.

\Proof. (a) The degree of the generator (D.1) is at least 
$2i_1 - 1 + 2 (j_1 + \ldots + j_s) > 2N - 1,$
and since it is odd, it is at least $2N+1$. On the other hand, it is at most 
$$ 1 + 3 + \ldots + (2N-1) + 2(j_1 + \ldots + j_s) 
\leq N^2 + 2N < (N+1)^2. $$

(b) If the degree is exactly $q = 2N+1$, then the argument in (a) implies 
$$ 2i_1 - 1 + 2 (j_1 + \ldots + j_s) = 2 N + 1, $$
which leads to $r = 1$ and $i_1 + j_1 + \ldots + j_s = N+1$. 
Only the partition 
$$ N+1 = (N+1) + 0 + \ldots +0 $$
does not contribute, and this proves (b). 

(c) If the degree $q$ of the generator in (D.1) is even and $> N$, then 
$j_1 + \ldots + j_s \leq N$ implies $r \geq 2$ and hence 
$$ q = (2i_1 - 1) + (2i_2 - 1) + 2 (j_1 + \ldots + j_s) 
> 2N-1 + 3 = 2N+2. 
\qeddis 

 Let $\tW$ be the graded vector space
$$ \tW = \bigoplus\limits_{q=2N}^{(N+1)^2-2} \tW_q ,$$
where $\dim \tW_q = \dim H^{q+1} (\WW)$. We view $\tW$ as a super vector space, where the parity
of $\tW_q$ is equal to the parity of $q$.

Let $\L_N$ be the free Lie superalgebra generated by $\tW$ (see e.g., [BMPZ92] for the
description of the basis of a free Lie superalgebra). The $\Z$-grading on $\tW$ extends to a 
$\Z$-grading of $\L_N$.

This Lie superalgebra was introduced by Haefliger [Hae76] in order to give the description
of the cohomology of the Lie algebra of smooth vector fields on a manifold with coefficients in the trivial 
module.

Let $\tilde V$ be a graded vector space obtained from $\L_N$ by a shift in the 
grading: $\tilde V_q = (\L_N)_{q-1}$.
Since the grading of $\L_N$ starts at component $2N$, the grading of $\tilde V$ starts at component $2N+1$.

 Let $\Omega^\bullet_M$ be the differential graded algebra of smooth 
differential forms on $M$, 
and consider the tensor product space $\Omega^\bullet_M \otimes \tilde V$ 
with the grading 
$$ \deg(\omega \otimes v) = \deg v - \deg \omega, \quad \omega \in \Omega^\bullet_M, v \in \tilde V.$$
 Since $\Omega^p_M$ vanishes for $p > N$, 
the grading on the space $\Omega^\bullet_M \otimes \tilde V$ starts at degree $N+1$. 
Next we consider the graded algebra $S^\bullet(\Omega^\bullet_M \otimes \tilde V)$ 
of supersymmetric multilinear forms on
 the graded superspace $\Omega^\bullet_M \otimes \tilde V$. We have
$$ S^\bullet(\Omega^\bullet_M \otimes \tilde V) = \mathop\oplus_{p=0}^\infty  
S^p(\Omega^\bullet_M \otimes \tilde V) ,$$
where $S^p(\Omega^\bullet_M \otimes \tilde V)_q = 0$ for $q<p \cdot (N+1)$. 
Hence, apart from the component of degree 0,
the grading on $S^\bullet(\Omega^\bullet_M \otimes \tilde V)$ starts in degree $N+1$.

\Theorem D.4. {\rm([Hae76], Theorem A)} 
Let $M$ be an $N$-dimensional smooth manifold. 
On the graded algebra $S^\bullet(\Omega^\bullet_M \otimes \tilde V)$ there exists a differential of degree 1, which
depends on a choice of the representatives of the Pontrjagin 
classes in $\Omega^\bullet_M$, and a homomorphism
of differential graded algebras
$$S^\bullet(\Omega^\bullet_M \otimes \tilde V) \rightarrow C^\bullet_c(\V_M,\R),$$
which induces an isomorphism in cohomology.
\qed

\Corollary D.5. $H^s_c(\V_M,\R) = 0$ \ for \ $1 \leq s \leq N$. 
\qed

If all Pontrjagin classes of $M$ vanish and there is a splitting algebra 
homomorphism $H^\bullet_M \into \Omega^\bullet_M$ 
(which is, for example, the case when $M$ is a torus or a compact Lie group), 
Haefliger gives a more
explicit realization for the cohomology of $\V_M$, which we described now. 

 Let $H^\bullet_M$ be the cohomology algebra 
of the manifold $M$. We view $H^\bullet_M$ as a commutative superalgebra   
and form the graded Lie superalgebra 
$$ \L (M) = H^\bullet (M,\R) \otimes \L_N $$
with the bracket of homogeneous elements given by 
$$ [\omega \otimes x, \omega^\prime \otimes x^\prime] 
= (-1)^{\deg x \cdot \deg \omega^\prime}
\omega \omega^\prime \otimes [x, x^\prime],$$
and the grading by 
$$ {\cal L}(M)_r = \sum_{p-q = r} H^p_M \otimes ({\cal L}_N)_q. $$ 
It follows in particular that $\L(M)_q = 0$ for $q<N$.
 
 Consider the cohomology $H^\bullet(\L(M))$ of the Lie superalgebra $\L(M)$ (see [Fu86], Sect.~1.6.3
for the definition of the cohomology for the Lie superalgebras, cf.\ Sect.~1.2 in [Hae76]).
Since $\L(M)$ is graded, its cohomology inherits a grading 
$H^p (\L(M),\R) = \oplus_q H^p (\L(M),\R)_q$ (see [Fu86], Sect.~1.3.7). 
Since the grading of $\L(M)$ starts at degree $N$,
we get that 
$$H^p (\L(M),\R)_q = 0 \quad {\rm for \ }  q < p\cdot N . \leqno{(D.2)}$$

\Theorem D.6. {\rm([Hae76], Theorem 3.4)} 
Let $M$ be a smooth manifold with a finite-dimensional
cohomology $H^\bullet_M$ for which all Pontryagin classes vanish. 
Suppose, moreover, that there exists a homomorphism
of differential graded algebras
$H^\bullet_M \rightarrow \Omega^\bullet_M,$
which induces an isomorphism in cohomology (this is the case if $M$ is a compact 
Lie group). Then
$$H_c^s ({\cal V}_M,\R) = \mathop\oplus\limits_{p+q=s} H^p (\L(M),\R)_q .
\qeddis 

\Corollary D.7. Under the assumptions of the preceding theorem, we have 
\litem{(a)} $H_c^s ({\cal V}_M,\R) = 0$ \ for \ $1 \leq s \leq N$.

\litem{(b)} For $N+1 \leq s \leq 2N+1$, $\dim H_c^s ({\cal V}_M,\R) 
= \sum\limits_{k=0}^N \dim H^{k}_M \cdot \dim \tW_{s+k-1}$. 

\litem{(c)} $\dim H_c^{N+1} ({\cal V}_M,\R) = p(N+1) - 1$ 
and only $H^N_M \cong \R$ contributes in 
the decomposition above. 
\litem{(d)} If $b_1(M) := \dim H^1_M = \dim H^{N-1}_M$, then 
$\dim H_c^{N+2} ({\cal V}_M,\R) = b_1(M) \cdot (p(N+1) - 1)$ and 
only $H^{N-1}_M \cong \R^{b_1(M)}$ contributes in the decomposition. 
 
\Proof. (a) If $H^p (\L(M),\R)_q$ is non-zero, then $q \geq p N$, so that (D.2) implies 
$s =p + q \geq (N+1)p \geq N+1$ if $s\not=0$. 

(b) We note that (D.2) implies $H_c^s ({\cal V}_M,\R) \cong H^1 (\L(M),\R)_{s-1}$ for $N+1 \leq s \leq 2N+1$. However,  
$$H^1 (\L(M),\R) = \left( \L(M) / [\L(M),\L(M)] \right)^* = \left( H^\bullet_M \otimes \tW \right)^* $$
and the claim (b) follows.
\qed

If $M$ 
is an $N$-dimensional torus $\T^N$, its cohomology $H^\bullet (\T^N)$ is an exterior algebra
in $N$ generators of degree 1. Combining the two previous corollaries, we get  
with $b_1(\T^N) = N$: 

\Corollary D.8. For the torus $T = \T^N$, we have: 
\litem{(a)} $H_c^s (\V_T,\R) = 0$ \ for \ $1 \leq s \leq N$. 
\litem{(b)} $\dim H_c^{N+1} (\V_T,\R) = p(N+1) - 1$. 
\litem{(c)} $\dim H_c^{N+2} (\V_T,\R) = N\cdot (p(N+1) - 1)$. 
\qed

\Remark D.9. For $M = \SS^n$, the algebra $H^\bullet({\cal V}_{\SS^n},\R)$ is 
calculated by Cohen and Taylor in [CT78]; Theorems 1.3 and 3.2. 
\qed

\def\entries{

\[AMR83 Abraham, R., J. E. Marsden, and T. Ratiu, ``Manifolds, Tensor Analysis, and Applications,'' 
Addison-Wesley, 1983 

\[ABFP05 Allison, B., S.~Berman, J.~Faulkner, and A.~Pianzola,
{\it Realizations of graded-simple algebras as loop algebras.}
math.RA/0511723.

\[BMPZ92 Bahturin, Y. A., Mikhalev, A. A., Petrogradsky, V. M., Zaicev, M. V.,
``Infinite-dimensional Lie superalgebras,''
Walter de Gruyter \& Co., Berlin, 1992

\[Beg87 Beggs, E. J., {\it The de Rham complex on infinite dimensional manifolds}, 
Quart. J. Math. Oxford (2) {\bf 38} (1987), 131--154 

\[BB99 Berman, S., and  Y.~Billig, 
{\it Irreducible representations for toroidal Lie algebras},  
J.\ Algebra {\bf 221} (1999), 188--231.

\[BR73 Bernshtein, I. N., and B. I. Rozenfel'd, {\it Homogeneous spaces 
of infinite-di\-men\-sio\-nal Lie algebras and characteristic classes of 
foliations}, Russ. Math. Surveys {\bf 28:4} (1973), 107--142 

\[BN06 Benkart, G., and E. Neher, {\it The centroid of extended affine 
and root graded Lie algebras}, J. Pure Appl. Algebra  {\bf 205:1} (2006), 117--145


\[Bi06 Billig, Y., 
{\it 
A category of modules for the full toroidal Lie algebra},  Int. Math. Res. Not. 2006, Art. ID 68395, 46 pp.

\[ChE48 Chevalley, C. and S. Eilenberg, {\it Cohomology theory of Lie groups and Lie 
algebras}, Transactions of the Amer. Math. Soc. {\bf 63} (1948), 85--124 
 
\[CT78 Cohen, F. R., and L. R. Taylor, {\it Computations of 
Gelfand-Fuks cohomology, the cohomology of function spaces, and the 
cohomology of configuration spaces}, in  ``Geometric applications 
of homotopy theory I'', Proc. Conf. Evanston, Ill., 106--173; 
Lectures Notes Math. {\bf 657}, Springer, 1978 

\[dWL83 de Wilde, M., and P. B. A. Lecomte, {\it Cohomology of the Lie algebra 
of smooth vector fields of a manifold, associated to the Lie derivative 
of smooth forms}, J. Math. pures et appl. {\bf 62} (1983), 197--214 

\[EM94 Eswara Rao, S., and Moody R.V.,
{\it Vertex representations for $n$-toroidal
Lie algebras and a generalization of the Virasoro algebra,}
Comm.Math.Phys. {\bf 159} (1994), 239-264.

\[FF01 Feigin, B. L. and D. B. Fuchs, {\it Cohomologies of Lie Groups and 
Lie Algebras}, in ``Lie Groups and Lie Algebras II'', A. L. Onishchik and 
E. B. Vinberg (Eds.), Encyclop. Math. Sci. {\bf 21}, Springer-Verlag, 2001 

\[FL80 Flato, M. and A. Lichnerowicz, {\it Cohomologie des repr\'esentations 
d\'efinies par la d\'erivation de Lie et \`a valeurs dans les formes, 
de l'alg\`ebre de Lie des champs de vecteurs d'une vari\'et\'e 
diff\'erentiable. Premiers espaces de cohomologie. Applications}, 
C. R. Acad. Sci. Paris S\'er. A-B {\bf 291:4}  (1980),  A331--A335

\[Fu86 Fuks, D. B. ``Cohomology of Infinite-Dimensional Lie Algebras,'' 
Consultants Bureau, New York, London, 1986 


\[GF70a Gelfand, I.\ M., Fuks, D.\ B., {\it Cohomology of the 
Lie algebra of vector fields with 
nontrivial coefficients}, Func.\ Anal. and its Appl. {\bf 4} (1970), 181--192

\[GF70b ---, {\it Cohomology of the Lie algebra of formal vector fields},
Izv. Akad. Nauk SSSR {\bf 34} (1970), 322--337


\[Go74 Godbillon, C., {\it Cohomologies d'alg\`ebres de Lie de champs de vecteurs formels},
S\'eminaire Bourbaki (1972/1973), Exp.No.421, Lecture Notes in Math. {\bf 383} (1974), 69--87



\[Hae76 Haefliger, A., {\it Sur la 
cohomologie de l'alg\`ebre de Lie des champs de vecteurs}, 
Ann. Sci. Ec. Norm. Sup. $4^e$ s\'erie {\bf 9} (1976), 503--532 


\[HS53] Hochschild, G., and J.-P. Serre, 
{\it Cohomology of Lie algebras}, Annals of Math. {\bf 57:3} (1953), 591--603 


\[Ka84 Kassel, C., 
{\it K\"ahler differentials and coverings of complex simple
Lie algebras extended over a commutative ring,} 
J.Pure Applied Algebra {\bf 34} (1984), 265--275 

\[Ko74 Koszul, J.-L., {\it Homologie des complexes de formes diff\'erentielles 
d'ordre \break sup\'erieur}, in ``Collection of articles dedicated to Henri Cartan 
on the occasion of his 70th birthday, I'', 
Ann. Sci. \'Ecole Norm. Sup. (4) {\bf 7} (1974), 139--153 


\[L99 T.A.~Larsson, 
{\it Lowest-energy representations of non-centrally extended diffeomorphism 
algebras.}
Commun.Math.Phys. {201}, (1999) 461-470.

\[Ma02 Maier, P., {\it Central extensions of topological
current algebras\/}, in ``Geometry and Analysis on Finite-
and Infinite-Dimensional Lie Groups,'' A. Strasburger et al Eds., 
Banach Center Publications {\bf 55}, Warszawa 2002; 61--76 



\[Ne04 Neeb, K.-H., {\it Abelian extensions of infinite-dimensional Lie groups}, 
Travaux math\'ematiques {\bf 15} (2004), 69--194 

\[Ne06a ---, {\it Lie algebra extensions and higher order cocycles}, 
J. Geom. Sym. Phys. {\bf 5} (2006), 48--74

\[Ne06b ---, {\it Non-abelian extensions of topological Lie algebras}, 
Communications in Algebra {\bf 34} (2006), 991--1041 



\[Neh04 Neher, E., {\it Extended affine Lie algebras}, 
C. R. Math. Acad. Sci. Soc. R. Can.  {\bf 26:3}  (2004), 90--96


\[PS86 Pressley, A., and G. Segal, ``Loop Groups,'' Oxford University Press, 
Oxford, 1986

\[Ro71 Rosenfeld, B. I., {\it Cohomology of certain infinite-dimensional 
Lie algebras}, Funct. Anal. Appl. {\bf 13} (1971), 340--342 

\[Tsu77 Tsujishita, T., {\it On the continuous cohomology of the 
Lie algebra of vector fields}, Proc. Jap. Math. Soc. {\bf 53:A} (1977), 
134--138 
 
\[Tsu81 ---, ``Continuous cohomology of the Lie algebra of vector fields,'' 
Memoirs of the Amer. Math. Soc. {\bf 34:253}, 1981
 
}

\references 

\def\address
{Yuly Billig 

School of Mathematics and Statistics

Carleton University

1125 Colonel By Drive

Ottawa, Ontario, K1S 5B6

Canada

billig@math.carleton.ca}

\def\addresstwo 
{Karl-Hermann Neeb

Technische Universit\"at Darmstadt 

Schlossgartenstrasse 7

D-64289 Darmstadt 

Deutschland

neeb@mathematik.tu-darmstadt.de}

\dlastpage 

\vfill\eject
\bye